\documentclass[a4paper]{amsart}
\usepackage[T1]{fontenc}
\usepackage{lmodern}
\usepackage[utf8]{inputenc}
\usepackage{amssymb}
\usepackage[all]{xy}
\usepackage{enumitem}
\usepackage{mathrsfs,dsfont,mathtools}
\usepackage{nicefrac}

\usepackage[pdftitle={Duality, correspondences and the Lefschetz 
map in equivariant KK-theory: a survey},
  pdfauthor={Heath Emerson and Ralf Meyer},
  pdfsubject={Mathematics; MSC 19K35, 46L80}
]{hyperref}
\newcommand*{\MRref}[2]{ \href{http://www.ams.org/mathscinet-getitem?mr=#1}{MR \textbf{#1}}}
\newcommand*{\arxiv}[1]{\href{http://www.arxiv.org/abs/#1}{arXiv: #1}}

\usepackage[lite]{amsrefs}
\usepackage{microtype}
\numberwithin{equation}{section}
\theoremstyle{plain}
\newtheorem{theorem}[equation]{Theorem}
\newtheorem{conjecture}[equation]{Conjecture}
\newtheorem{lemma}[equation]{Lemma}
\newtheorem{proposition}[equation]{Proposition}
\newtheorem{corollary}[equation]{Corollary}
\theoremstyle{definition}
\newtheorem{definition}[equation]{Definition}
\newtheorem{notation}[equation]{Notation}
\theoremstyle{remark}
\newtheorem{remark}[equation]{Remark}
\newtheorem{example}[equation]{Example}

\DeclareMathOperator{\Hom}{Hom}
\DeclareMathOperator{\KK}{KK}
\DeclareMathOperator{\RKK}{RKK}

\DeclareMathOperator{\K}{K}

\DeclareMathOperator{\RK}{RK}
\DeclareMathOperator{\vbK}{VK}
\DeclareMathOperator{\Vect}{Vect}

\newcommand*{\GKK}{\widehat{\textsc{kk}}{}}

\newcommand*{\ev}{\textup{ev}}

\newcommand*{\PD}{\textup{PD}}

\newcommand*{\forget}{\textup{forget}}

\newcommand*{\op}{\textup{op}}

\newcommand*{\Cliff}{\textup{Cliff}}

\newcommand*{\Q}{\mathbb{Q}}

\newcommand*{\C}{\mathbb C}
\newcommand*{\Z}{\mathbb Z}

\newcommand*{\R}{\mathbb R}
\newcommand*{\Comp}{\mathbb K}

\newcommand*{\Index}{\mathscr{I}}

\newcommand*{\Hils}{\mathcal{H}}
\newcommand*{\CONT}{\mathcal{C}}
\newcommand*{\EG}{\mathcal E}

\newcommand*{\Tvert}{\textup T}
\newcommand*{\tang}{\Tvert \Tot} 
\newcommand*{\Grd}{\mathcal G}
\newcommand*{\ind}{\mathrm{Ind}}
\newcommand*{\Base}{Z} 
\newcommand*{\Tot}{X}  
\newcommand*{\Other}{Y}
\newcommand*{\Source}{X}
\newcommand*{\Target}{Y}
\newcommand*{\Third}{U}
\newcommand*{\Dual}{P}
\newcommand*{\Midd}{M}
\newcommand*{\VB}{V}
\newcommand*{\Triv}{E}
\newcommand*{\Normal}{\textup N}
\newcommand*{\Kclass}{\xi}
\newcommand*{\ctau}{\mathscr{C}_\tau}

\newcommand*{\NM}{\Phi}
\newcommand*{\anchor}{\varrho}
\newcommand*{\mapl}{b}
\newcommand*{\mapr}{f}
\newcommand*{\pr}{\mathrm{pr}}

\newcommand*{\UNIT}{\mathds{1}}     

\newcommand*{\nb}{\nobreakdash}
\newcommand*{\Cst}{\textup{C}^*}
\newcommand*{\Star}{$^*$\nobreakdash-\hspace{0pt}}

\newcommand*{\total}[1]{\lvert#1\rvert}
\newcommand*{\proj}[1]{\pi_{#1}}
\newcommand*{\zers}[1]{\zeta_{#1}}
\newcommand*{\norm}[1]{\lVert#1\rVert}

\newcommand*{\defeq}{\mathrel{\vcentcolon=}}

\newcommand*{\comul}{\nabla}
\newcommand*{\blank}{\textup{\textvisiblespace}}

\newcommand*{\opem}{\hookrightarrow}

\newcommand*{\Stab}{\textup{Stab}}
\newcommand*{\selfmap}{\varphi}
\newcommand*{\Selfmap}{\varphi}
\newcommand*{\mult}{\mathbf{n}}
\newcommand*{\Fix}{\textup{Fix}}
\newcommand*{\sign}{\textup{sign}}
\newcommand*{\ID}{\textup{id}}

\newcommand*{\Rep}{\textup{Rep}}

\newcommand*{\Fixed}{\textup{F}}

\DeclareMathOperator{\Lef}{Lef}
\DeclareMathOperator{\Eul}{Eul}

\begin{document}

\title[Duality and correspondences]{Duality, correspondences and the Lefschetz map in equivariant KK-theory: a survey}

\author{Heath Emerson}
\email{hemerson@math.uvic.ca}

\address{Department of Mathematics and Statistics\\
  University of Victoria\\
  PO BOX 3045 STN CSC\\
  Victoria, B.C.\\
  Canada V8W 3P4}

\begin{abstract}
We survey work by the author and Ralf Meyer on equivariant 
KK-theory. Duality plays a key role in our approach. 
We have organized this survey around the objective of computing 
a certain homotopy-invariant of a 
space equipped with a (generally proper) action of a groupoid. This invariant is  
called the Lefschetz map. The Lefschetz map associates an equivariant K-homology 
class to an equivariant Kasparov self-morphism of a space X.
We want to describe it explicitly in the setting of bundles of smooth manifolds 
over the base space of a groupoid, in which groupoid elements act by 
diffeomorphisms between fibres. 
To get the required description we develop a topological model of equivariant 
KK-theory by way of a theory of correspondences, building on ideas 
of Paul Baum, Alain Connes and Georges Skandalis in the 1980's. 
This model agrees with the analytic model for bundles of smooth manifolds 
under some technical conditions 
related to the existence of equivariant vector bundles. Subject 
to these conditions we obtain the desired computation of the Lefschetz map  
in purely topological terms. Finally, we describe a generalization of the 
classical Lefschetz fixed-point formula to apply to correspondences, 
instead of just maps. 

 \end{abstract}
\subjclass[2000]{19K35, 46L80}
\thanks{Heath Emerson 
  was supported by a National Science and Engineering Council 
  of Canada (NSERC) Discovery grant. }
\maketitle

The papers \cite{Emerson-Meyer:Euler},  \cite{Emerson-Meyer:Dualities},
\cite{Emerson-Meyer:Equivariant_K},  \cite{Emerson-Meyer:Geometric_KK} 
present a study of the equivariant Kasparov groups \(\KK^\Grd\bigl( \CONT_0(X), \CONT_0(Y)\bigr)\)
where \(\Grd\) is a 
locally compact Hausdorff groupoid with Haar system and \(X\) and 
\(Y\) are \(\Grd\)-spaces, usually with \(\Tot\) a proper \(\Grd\)-space. 
 This program builds on work of Kasparov, Connes and 
Skandalis done mainly in the 1980's. At that point, the main interest was the index theorem
of Atiyah and Singer and its generalisations, and later, the Dirac dual-Dirac method 
and the Novikov conjecture. For us, the goal is to develop
Euler characteristics and Lefschetz formulas in equivariant \(\KK\)-theory. Via the 
Baum-Connes isomorphism -- when it applies -- this contributes to noncommutative 
topology and index theory. Our program started in 
\cite{Emerson-Meyer:Euler} where we found the \emph{Lefschetz map} in connection 
with a \(\K\)-theory problem. 
 We will give the 
definition of the Lefschetz map 
in the first section, but for now record that it has the form  
\begin{equation}
\label{eq:lefmap}
\Lef \colon  \KK^{\Grd\ltimes \Tot}_*\bigl(  \CONT_0(\Tot\times_\Base \Tot), \CONT_0(\Tot)\bigr) \to 
\KK_*^\Grd(\CONT_0(\Tot), \CONT_0(\Base)\bigr),
\end{equation}
where we always denote by \(\Base\) the base space of the groupoid.  This 
map is defined under certain somewhat technical circumstances, but, again, these 
normally involve proper \(\Grd\)-spaces \(X\).
The domain of the Lefschetz map is very closely related to the simpler-looking group 
\(\KK_*^\Grd\bigl( \CONT_0(\Tot), \CONT_0(\Tot)\bigr)\): the  latter group maps 
in a natural way to the domain in \eqref{eq:lefmap} and this map is an isomorphism 
when the anchor map \(X \to \Base\) is a proper map. This means that the Lefschetz 
map can be used to assign an invariant, which is an equivariant \(\K\)-homology class, to an
 equivariant Kasparov self-morphism of \(\Tot\). We call this class the 
 \emph{Lefschetz invariant} of the map. It  bears consideration even when 
 \(\Grd\) is the trivial groupoid, and the reader can do worse than to consider this 
 case to begin with, although by doing so one misses the applications to 
 noncommutative topology.

 
 The definition of \eqref{eq:lefmap} uses the notion of an   
 \emph{abstract dual} for \(\Tot\). 
 Abstract duals for a given \(\Grd\)-space \(X\) are not unique but the Lefschetz 
 map does not depend on the choice of a dual, only on the existence of one.
 Abstract duals do not always exist either: a Cantor set \(\Tot\) doesn't have an 
 abstract dual even if \(\Grd\) is trivial.  But if \(X\) is a smooth \(\Grd\)-manifold, 
 with \(\Grd\) acting smoothly and properly on \(X\), then \(X\) has an 
 abstract dual, and if \(\Grd\) is a group, then any \(\Grd\)-simplicial complex has an abstract dual
 due to \cite{Emerson-Meyer:Euler}.

The
 Lefschetz map is functorial for \(\Grd\)-maps \(\Tot \to \Tot'\) in a way 
 made explicit in \S \ref{subsec:the_lefschetz_map} (see  
 Theorem \ref{thm:Lef_well-defined}). In brief, 
it is a \emph{homotopy invariant} of the \(\Grd\)-space \(X\). 
  In particular, since 
 \(\Lef\) doesn't depend on the dual used to
  to compute it,  one can try to compute the Lefschetz invariant of 
  a given morphism using two different duals and thereby get 
  an identity in equivariant 
  \(\K\)-homology (see \S \ref{subsec:map_examples}).
 Such examples (worked out in \cite{Emerson-Meyer:Euler} and 
 \cite{Emerson-Meyer:Equi_Lefschetz}), 
 seemed to us interesting enough to support making a systematic 
 study of the Lefschetz map. However, to get started on this question one obviously 
 has to first describe the morphisms
 \(f\) themselves in some kind of satisfactory way.  To this end, we have 
 extended the theory of correspondences initiated by Baum, Connes, 
 Skandalis and others, to the equivariant situation, in the paper 
 \cite{Emerson-Meyer:Geometric_KK}. This extension 
 presents some new features, and we will devote a part of this survey 
 to explaining them. Of course the theory of correspondences is 
 useful and important in its own right. But it is \emph{designed} for 
 intersection theory because of the way correspondences are composed 
 using coincidence spaces and transversality. 
 
The theory of correspondences requires the groupoid \(\Grd\) to be \emph{proper}. 
The Baum-Connes conjecture allows us to reduce to this case subject to 
a weaker assumption that we explain below. Let \(\Grd\) be proper and 
\(X\) and \(Y\) be \(\Grd\)-spaces. A
 \emph{\(\Grd\)-equivariant correspondence from \(X\) to \(Y\)} is a 
quadruple \((M, \mapl, \mapr, \xi)\) where \(M\) is a \(\Grd\)-space, 
\(\mapl\colon M \to X\) is a 
\(\Grd\)-map (not necessarily proper), \(\xi\) is an equivariant \(\K\)-theory class
with compact vertical support along the fibres of \(b\), and 
\(f\) is a \(\K\)-oriented 
\emph{normally non-singular map } from \(M\) to \(Y\) (see Definition \ref{def:normal_map}). 
For example, if
 \(\Grd\) is a compact group,  \(Y\) is a point and \(M\) is compact, then 
an normally non-singular map  \(M\to Y\) 
is the specification of an orthogonal representation of \(\Grd\) on some \(\R^n\),
an equivariant vector bundle \(V\) over \(M\), and an open equivariant embedding 
\(\hat{f}\colon V \to \R^n\).  To construct an example of 
such a triple, assume that \(M\) has been given the structure of a smooth manifold, 
and that \(\Grd\) acts smoothly. In this case we may appeal to a theorem of Mostow to 
embed \(M\) in a finite-dimensional 
linear representation of \(\Grd\), then take \(\VB\) to be the normal bundle to 
the embedding.

There is a topologically defined equivalence relation on correspondences that makes the 
set of equivalence classes of \(\Grd\)-equivariant correspondences from 
\(X\) to \(Y\) the morphism set \(\GKK_\Grd^*(X,Y)\) 
in a \(\Z/2\)-graded category \(\GKK^\Grd\) which maps naturally to \(\KK^\Grd\).
For example, let \(\Grd\) be  compact group and let 
both \(X\) and \(Y\) be the one-point space.
Let \(M\) be a smooth, compact, equivariantly \(\K\)-oriented, even-dimensional
 \(\Grd\)-manifold, 
\(\xi \in \K^0_\Grd (M)\) be an equivariant \(\K\)-theory class for \(M\) 
represented by an 
equivariant vector bundle \(\VB\) on \(M\). By embedding  \(M\) in a finite-dimensional 
representation of \(\Grd\) as in the previous paragraph, we can endow the map 
from \(M\) to \(Y\defeq \star\) with the structure of a smooth, \(\K\)-oriented, normal 
map, and we obtain a \(\Grd\)-equivariant correspondence 
\((M, \star, \star, \xi) \) from a point to itself. This  yields a class in 
\(\GKK_\Grd^0( \star, \star)\). Applying the natural map 
\(\GKK^\Grd_0(\star, \star) \to 
\KK^\Grd_0(\C, \C) \cong \Rep (\Grd)\) maps this correspondence to 
the \emph{\(\Grd\)-equivariant topological index of \(D_V\)} in the sense of 
\cite{Atiyah-Singer:Index_I}, where \(D_V\) is the Dirac operator on \(M\) 
twisted by the equivariant vector bundle \(\VB\). By the Atiyah-Singer
Index theoreom, this agrees with 
the \(\Grd\)-equivariant analytic index of \(D_\VB\) in 
\(\Rep (\Grd)\), obtained by considering the difference of 
finite-dimensional \(\Grd\)-representations on the kernel and cokernel of 
\(D_\VB\). 

The combination of the Atiyah-Singer index theorem and the theory of equivariant 
correspondences represents a powerful tool, because while the index theorem 
allows us to translate analytic problems into topological ones, the theory of 
correspondences allows us to manipulate this topological data in interesting 
ways. For an example of this process in connection with the
 representation theory of complex semisimple 
Lie groups, see \cite{Emerson-Yuncken:BBW}. 

In terms of the Lefschetz map, the fact that 
correspondences can be composed in an essentially topological fashion
has the consequence that the Lefshetz invariants of self-correspondences
of a smooth \(\Grd\)-manifold \(\Tot\), or, 
or more precisely, of their images in \(\KK^\Grd\), can be computed in 
terms of considerations of transversality. 
We explain the outcome of this computation in 
Section \ref{sec:dualities_and_lef}: the gist is  
that the Lefschetz invariant of a 
smooth equivariant self-correspondence \(\Psi\)
of a smooth \(\Grd\)-manifold \(X\) in general position, 
can be described in terms of a \(\Grd\)-space called
the \emph{coincidence space} \(\Fixed_\Psi'\) of the 
correspondence. The coincidence space
 inherits from the smooth structure and \(\K\)-orientation on 
\(\Psi\) the structure of a 
 smooth and equivariantly \(\K\)-oriented \(\Grd\)-manifold  
 which maps to \(\Tot\) and represents
an equivariant correspondence from \(\Tot\) to \(\Base\), thus a cycle for
\(\GKK_\Grd^*(\Tot, \Base)\) and then a class in
 \(\KK_*^\Grd\bigl( \CONT_0(\Tot), \CONT_0(\Base)\bigr)\). It represents 
 the Lefschetz invariant of the morphism represented by \(\Psi\). 
  See Theorem 
 \ref{thm:main_lef} for the exact statement. The `general position' caveat
 is non-trivial: in the equivariant setting, it may not be possible to 
 perturb a pair of equivariant maps to make them transverse. The 
 framework of equivariant correspondences allows us to treat this 
 difficulty using Bott Periodicity, but we do not discuss this much here.

 For example, if \(\Grd\) is a 
 compact group, \(\Tot\) a smooth and compact manifold with a smooth action of 
 \(\Grd\), then the Lefschetz invariant of a smooth equivariant map 
 \(f\colon X\to X\) in general position is the fixed-point set of the map, which is a 
 finite set of points permuted by \(\Grd\), oriented by an equivariant 
 line bundle over this finite set. This bundle depends on orientation data 
 from the original map \(f\) in a manner which reduces to the classical 
 choice of signs at each fixed-point when \(\Grd\) is trivial. Thus, 
 the topological model of the Lefschetz map provided by the theory of 
 correspondences yields an interpretation of \(\Lef\) in terms of a 
 a fixed-point theory for correspondences.

One naturally asks when the map 
\(\GKK_\Grd^*(X, Y)\to \KK_*^\Grd\bigl( \CONT_0(X), \CONT_0(Y)\bigr)\) is an 
isomorphism.  We explain our results on this in \S \ref{subsec:normal_finite__type};
 once again, they 
rely on duality in a crucial way. When they apply, the topological and analytic 
Lefschetz maps are equivalent. As mentioned above, 
the Baum-Connes conjecture can be used to reduce the non-proper situation 
to the proper one under some weaker assumptions on the \(\Grd\) action on 
\(X\), namely that it be \emph{topologically amenable}.
 This is explained in \S \ref{subsec:inflation_trick}. Putting everything together 
gives a computation of the Lefschetz invariant for quite a wide spectrum of 
smooth \(\Grd\)-spaces \(X\). 

What is duality? It is central to our whole framework, and accordingly 
we begin the article with a 
discussion of it.  
It is well-known from the work of Kasparov and Connes and
Skandalis (see \cite{Kasparov:Novikov} and \cite{Connes-Skandalis:Longitudinal})
that if if \(\Tot\) is a smooth 
manifold, then there
is a natural family of isomorphisms \[ \KK(\CONT_0(\tang),\C)
 \cong \RK^*(X) \defeq \KK^\Tot (\CONT_0(X), \CONT_0(X) )\]
where the groupoid equivariant \(\KK\) group on the right is equivariant 
\emph{representible} \(\K\)-theory, or \(\K\)-theory with \emph{locally finite support}, 
denoted \(\RK^* (X)\) by Kasparov.
 There is a generalisation 
of this duality to the equivariant situation if \(\Grd\) is a groupoid acting 
smoothly and properly on a \emph{bundle} \(X\to \Base\) of smooth 
manifolds over the base \(\Base\) of \(\Grd\),
 and furthermore, the roles of \(X\) and \(\tang\) can be 
in a sense reversed, so that one can establish a pair of natural (in a 
technical sense) familes of 
isomorphisms
\begin{equation}
\label{eq:first_duality}
 \KK_*^\Grd( \CONT_0(\tang)\otimes A, B) \cong 
\KK^{\Grd\ltimes X}_* ( \CONT_0(X)\otimes A, \CONT_0(X)\otimes B)
\end{equation}
and 
\begin{equation}
\label{eq:second_duality}
 \KK^{\Grd}_*(\CONT_0(X)\otimes A, B) 
\cong \KK^{\Grd\ltimes X}_*(\CONT_0(X) , \CONT_0(\tang)\otimes B)
\end{equation}
for all \(\Grd\)-C*-algebras \(A\) and \(B\).  (The tensor products are all 
in the category of \(\Grd\)-C*-algebras.) 
 These results are proved in \cite{Emerson-Meyer:Dualities}. 
 It is the first kind of duality \eqref{eq:first_duality} which is relevant for the 
 Lefschetz map, and the second \eqref{eq:second_duality} that is used 
 to prove that the map 
 \(\GKK^\Grd\to \KK^\Grd\) is an isomorphism in certain cases. 
 The basic idea is that since 
 the duality isomorphisms \eqref{eq:second_duality} are 
 themselves induced by equivariant correspondences, duality 
 can be used simultaneously in both the analytic and 
 topological categories to to reduce 
 the question to a problem about \emph{monovariant} \(\KK\)-theory, that is, 
equivariant \(\K\)-theory with support conditions. 

What is needed to make this work is then a topological model of 
duality. The main new issue that appears is that our equivariant 
correspondences require a good supply of equivariant vector bundles 
and this forces conditions on the groupoid \(\Grd\). These  
considerations have in fact already appeared in the literature in 
connection with the (proper) groupoids \(G\ltimes \EG{G}\) in work by  
Wolfang L\"uck and Bob Oliver in \cite{Lueck-Oliver:Completion} (see 
\S \ref{subsec:normal_maps} for the details) where \(G\) is a 
discrete group and \(\EG{G}\) is its classifying space for proper 
actions. We explain exactly what the conditions are and how they 
are related to embedding theorems generalizing the embedding theorem 
of Mostow alluded to above. 

The classical Lefschetz fixed-point theorem relates fixed-points of a map 
and the homological invariant of the map 
obtained by taking the graded trace of 
the induced map on homology, called the Lefschetz number. 
We finish this survey by introducing some global, homological invariants of 
correspondences which generalize the Lefschetz number, at least in the 
case when the groupoid \(\Grd\) is trivial. 
Roughly speaking, a
Kasparov self-morphism, and in particular a self-correspondence
should be considered as determining a linear 
\emph{operator} on homology 
instead of just a number. We call it the Lefschetz operator. We will explain how to  
extend the classical Lefschetz fixed-point theorem to correspondences by 
identifying the Lefschetz operator with the operator of pairing with the 
Lefschetz invariant. In situations where a local index formula is available, 
this results in a description of the Lefschetz operator in local, geometric 
terms. 

It remains to describe the Lefschetz map in global, homological terms -- as 
in the classical formula, in which fixed-points are related to traces on 
homology. This problem seems quite delicate, however. One way of 
proceeding is to replace \(\K\)-theory groups by \(\Grd\)-equivariant \(\K\)-theory 
\emph{modules} over \(\Rep (\Grd)\) and replace the ordinary trace by the 
Hattori-Stallings trace. However, this is only defined under quite stringent 
conditions: such modules are only finitely presented in general only for 
groups for which the representation rings have finite cohomological 
dimension, and this requires additional hypotheses on \(\Grd\). This 
work is still in progress. 

\section{Abstract duality and the Lefschetz map}
\label{sec:abstract_duals}

Throughout this paper, \emph{groupoid} shall mean locally compact 
Hausdorff groupoid with Haar system. All topological spaces will be 
assumed paracompact, locally compact and Hausdorff.  For the material in this section, 
see \cite{Emerson-Meyer:Dualities}. For source material on equivariant 
\(\KK\)-theory for groupoids, see \cite{LeGall:KK_groupoid}. 
One seems to be forced to consider groupoids, as opposed to groups, in 
equivariant Kasparov theory, even if one is ultimately only interested in 
groups. This will be explained later. 

Therefore we will work more or less uniformly 
with groupoids when discussing general theory. When we discuss 
topological equivariant Kasparov theory, we will further assume that 
all groupoids are \emph{proper}. This restriction is needed for various 
geometric constructions. The additional assumption of properness 
involves no serious loss of generality for our purposes 
because the Baum-Connes 
isomorphism, when it applies, gives a method of replacing non-proper 
groupoids by proper ones.

\subsection{Equivariant Kasparov theory for groupoids}
Let \(\Grd\) be a groupoid. We let \(\Base\) denote the base space.
A \(\Grd\)-C*-algebra is in particular 
a C*-algebra over \(\Base\). This means that there is given a non-degenerate 
equivariant *-homomorphism 
from \(\CONT_0(\Base)\) to central multipliers of \(A\).
This identifies \(A\) with the section algebra of a continuous 
bundle of C*-algebras over \(\Base\). For a groupoid action we require 
in addition an isomorphism 
\( r^*(A) \to s^*(A)\) which is compatible with the structure of 
\(r^*(A)\) and \(s^*(A)\) as C*-algebras over \(\Grd\). 
Here \(r\colon \Grd\to \Base\) and \(s\colon \Grd\to \Base\) are the range and 
source map of the groupoid, and 
\(r^*\) (and similarly \(s^*\)) denotes the usual pullback operation of bundles. 
From the bundle point of view, all of this means that 
groupoid elements \(g\) with \(s(g) = x\) and \(r(g) = y\) induce 
*-homomorphisms \(A_x \to A_y\) between the fibres of 
\(A\) at \(x\) and \(y\). 

In particular, if \(A\) is commutative, then \(A\) is the C*-algebra 
of continuous functions on a locally compact \(\Grd\)-space \(X\),
equipped with a map \(\anchor_X\colon X \to \Base\) called the 
\emph{anchor map} for \(X\), and a homeomorphism 
\[ \Grd\times_{\Base,s} X  \to \Grd\times_{\Base,r} X, \hspace{1cm} (g,x) \mapsto (g, gx)\]
where the domain and range of this homeomorphism (by abuse of notation) 
are respectively 
\[ \Grd\times_{\Base, s} X \defeq \{ (g, x) \in \Grd\times X \mid s(g) = \anchor_X (x)\},\]
and similarly for \(\Grd\times_{\Base, r} X\) using \(r\) instead of \(s\).

\subsection{Tensor products}
The category of \(\Grd\)-C*-algebras has a symmetric monoidal structure
given by tensor products. We describe this very 
briefly (see \cite{Emerson-Meyer:Dualities}[Section 2] for 
details).

Let \(A\) and \(B\) be two \(\Grd\)-C*-algebras. Since they are each 
C*-algebras over \(\Base\), their external tensor product \(A\otimes B\) 
is a C*-algebra over \(\Base\times \Base\). We restrict this to 
a C*-algebra over the diagonal \(\Base\subset \Base\times \Base\).
The result is called the \emph{tensor product of \(A\) and \(B\) over \(\Base\)}. 
The tensor product of \(A\) and \(B\) over \(\Base\) carries a diagonal action of 
\(\Grd\). We leave it to the reader to check that we obtain a \(\Grd\)-C*-algebra
in this way. 

In order not to complicate notation, \emph{we write just 
\(A\otimes B\) for the tensor product of \(A\) and \(B\)
 in the category of \(\Grd\)-C*-algebras.} We emphasize that the 
 tensor product is over \(\Base\); this is not the same as the tensor 
 product in the category of C*-algebras. 
 
 For commutative C*-algebras, \emph{i.e.} for \(\Grd\)-spaces, say 
 \(X\) and \(Y\), with anchor maps as usual denoted 
 \(\anchor_X\colon X \to \Base\) and \(\anchor_Y\colon Y \to \Base\), 
 the tensor product is Gelfand dual to the operation which forms from \(X\) and 
 \(Y\) the fibre product  
 \[ X\times_\Base Y \defeq \{ (x,y)\in X\times Y \mid \anchor_X(x) = \anchor_Y (y)\}.\]
 The required 
 anchor map \(\anchor_{X\times_\Base Y}\colon X \times_\Base Y \to \Base\) 
 is of course the composition of the first coordinate projection and the 
 anchor map for \(X\) (or the analogue using \(Y\); they are equal). 
 Of course groupoid elements act diagonally in 
 the obvious way. Such \emph{coincidence spaces} as the one just 
 described will appear again and again in the theory of correspondences.

Finally, for the record, we supply the following important 
definition.

\begin{definition}
\label{def:proper_action}
Let \(\Grd\) be a groupoid. A \(\Grd\)-space \(X\) is \emph{proper} if the map 
\[ \Grd\times_{\Base} X\to X, \; (g, x) \mapsto (gx, x) \]
is a proper map, where \(\Grd\times_\Base X \defeq \{ (g, x)\in \Grd\times X\mid 
s(g) = \anchor_X (x)\}\). 
\end{definition}

A groupoid is itself called \emph{proper} if it acts properly on its base space 
\(\Base\). Explicitly, the map 
\[ \Grd \to X\times_\Base X, \hspace{1cm} g \mapsto (r(g), s(g))\]
is required to be proper.

\subsection{Equivariant Kasparov theory}
Le Gall has defined \(\Grd\)-equivariant \(\KK\)-theory in 
\cite{LeGall:KK_groupoid}. We briefly sketch the definitions. 
Let 
\(A\) and \(B\) be (possibly \(\Z/2\)-graded) 
\(\Grd\)-C*-algebras. 

Then a cycle for 
\(\KK^\Grd (A,B) \) is given by a \(\Z/2\)-graded 
\(\Grd\)-equivariant Hilbert \(B\)-module 
\(\mathcal{E}\), together with a \(\Grd\)-equivariant 
grading-preserving 
*-homomorphism from \(A\) to the C*-algebra of bounded, 
adjointable operators on \(\mathcal{E}\), and 
an essentially \(\Grd\)-equivariant self-adjoint 
operator \(F\) on \(\mathcal{E}\) which is graded odd
 and satisfies 
\( [a,F]\) and \(a (F^2-1)\) are compact operators (essentially 
zero operators) for all \(a\in A\). 

Modulo an appropriate equivalence 
relation, the set of equivalence classes of cycles 
can be identified with the morphism set 
\(\KK^\Grd(A,B)\) in an additive, symmetric monoidal 
category.  Higher \(\KK\)-groups are defined using 
Clifford algebras, and since these are \(2\)-periodic, 
there are only two up to isomorphism. 
 We denote by 
\(\KK_*^\Grd(A,B)\) the sum of these two groups.  

If \(A\) and \(B\) are \(\Grd\)-C*-algebras, then 
the group \(\RKK^\Grd(\Tot; A
,B)\) is by definition the groupoid-equivariant Kasparov group 
 \(\KK^{\Grd\ltimes \Tot}
( \CONT_0(X)\otimes A, \CONT_0(X)\otimes B).\)
The tensor products are in the category of \(\Grd\)-C*-algebras. 
This group differs from \(\KK^\Grd(A, \CONT_0(X)\otimes B)\) only 
in the support condition on cycles. For example if \(\Grd\) is trivial and 
\(A=B=\C\) then \(\KK^\Grd(A,\CONT_0(X)\otimes B)\) is the ordinary \(\K\)-theory 
of \(X\) and \(\RKK^\Grd(X; A,B)\) is the representable \(\K\)-theory of \(X\) (a non-compactly
supported theory.) We discuss these groups in more detail in the next section. 
Of course similar remarks hold for higher \(\RKK^\Grd\)-groups.

\subsection{Equivariant \(\K\)-theory}
\label{subsec:equivariant_K-theory}
In this section, we present an exceedingly brief overview of equivariant \(\K\)-theory, 
roughly sufficient for the theory of equivariant correspondences. For more details see
\cite{Emerson-Meyer:Equivariant_K}.

Let \(X\) be a proper \(\Grd\)-space.  Recall that a \(\Grd\ltimes X\)-space 
consists of a \(\Grd\)-space \(Y\) together with a \(\Grd\)-equivariant 
map \(\anchor_Y: Y \to X\) serving as the anchor map for the 
\(\Grd\ltimes X\)-action. 

\begin{definition}
\label{def:representible_K}
Let \(Y\) be a \(\Grd\ltimes X\)-space. The 
\emph{\(\Grd\)-equivariant representable \(\K\)-theory of \(Y\) with 
\(X\)-compact supports} is the group
\[\RK^{-*}_{\Grd, X} (Y) \defeq \KK^{\Grd\ltimes X}_* \bigl( \CONT_0(X), \CONT_0(Y)\bigr).\]

The \emph{\(\Grd\)-equivariant representable \(\K\)-theory of \(Y\)}
 is \[\RK_\Grd^{*} (Y) \defeq \RK^*_{\Grd, Y} (Y).\]

\end{definition}

Cycles for \(\KK^{\Grd\ltimes X} \bigl( \CONT_0(X), \CONT_0(Y)\bigr)\)
consist of pairs \( (\Hils, F)\) where 
\(\mathcal{E}\) is a countably generated \(\Z/2\)-graded 
\(\Grd\ltimes X\)-equivariant right Hilbert \(\CONT_0(X)\)-module equipped with 
a \(\Grd\ltimes X\)-equivariant non-degenerate 
*-homomorphism from \(\CONT_0(X)\) to 
the C*-algebra of bounded, adjointable operators on \(\mathcal{E}\), 
and \(F\) is a bounded, odd, self-adjoint, 
essentially \(\Grd\)-equivariant
 adjointable operator on \(\Hils\) such that  
\( f (F^2-1)\) is a compact operator, for all 
\(f\in \CONT_0(X)\). The 
properness of \(\Grd\) implies that \(F\) may be averaged to be 
actually \(\Grd\)-equivariant, so we assume this in the following. 

The Hilbert \(\CONT_0(Y)\)-module \(\mathcal{E}\) is the space of 
continuous sections of a continuous 
field of \(\Z/2\)-graded Hilbert spaces \(\{ \Hils_y\mid y\in Y\}\) 
over \(Y\). Since \(F\) must be \(\CONT_0(Y)\)-linear, 
it consists of a continuous family \(\{F_y\mid y\in Y\}\) of odd operators on these 
graded Hilbert spaces such that \(F^2_y-1\) is a compact operator on \(\Hils_y\) 
for all \(y \in Y\). 

By \(\Grd\ltimes X\)-equivariance, the representation of 
\(\CONT_0(X)\) on \(\mathcal{E}\) must factor through the 
*-homomorphism \(\CONT_0(X) \to \CONT_0(Y)\) Gelfand dual to the 
anchor map \(\anchor_Y\colon Y \to X\). Therefore
 \(F\) commutes with the action of \(\CONT_0(X)\) as well; in fact 
 the induced representation of 
\(\CONT_0(X)\) on each  
Hilbert space \(\Hils_y\) sends a continuous function \(f\in \CONT_0(X)\) 
to the operator of multiplication by the complex number 
\(f\bigl( \anchor_Y (y)\bigr)\). 

In particular, the only role of the representation of 
\(\CONT_0(X)\) is to relax the support condition on the compact-operator 
valued-function \(F^2-1\) from requiring it to vanish at \(\infty\) of \(Y\) 
to only requiring it to vanish at infinity along the fibres of 
\(\anchor_Y\colon Y \to X\).

If \(\anchor_Y\colon Y \to X\) is a \emph{proper} map then 
\(\RK_{\Grd, X}^* (Y) = \RK^*_\Grd (X) \defeq \RK^*_{\Grd, X} (X)\); these two groups have 
exactly the same cycles.

\begin{example}
\label{ex:vector_bundles}
Any \(\Grd\)-equivariant complex vector bundle \(\VB\) 
on \(Y\) yields a cycle for \(\RK^0_\Grd (Y)\) by choosing a \(\Grd\)-invariant 
Hermitian metric on \(\VB\) and forming the corresponding 
\(\Grd\ltimes Y\)-equivariant \(\Z/2\)-graded 
right Hilbert \(\CONT_0(Y)\) module of sections, where the grading is 
the trivial one.  We set 
the operator equal to zero.
\end{example}

\begin{example}
\label{ex:Thom_isomorphism}
Let \(X\) be a \(\Grd\)-space and let \(\VB\) be a \(\Grd\)-equivariantly \(\K\)-oriented
 vector bundle over \(X\) of (real) dimension \(n\). 
 The \(\Grd\)-equivariant vector bundle projection \(\proj{\VB}\colon \VB \to X\) 
 gives \(\VB\) the structure of a space over \(X\), so that \(\VB\) becomes a 
 \(\Grd\ltimes X\)-space.  
 Then the \emph{Thom isomorphism} provides an invertible \emph{Thom class}
 \[ t_\VB \in \RK^{\dim \VB}_{\Grd, X} (\VB)\defeq \KK^{\Grd\ltimes X}_{\dim \VB}\bigl(\CONT_0(X), \CONT_0(\VB)\bigr).\]
In the case \(\Grd = \mathrm{Spin}^c(\R^n)\) and \(X=\star\) and \(\VB\defeq \R^n\) with 
the representation \(\mathrm{Spin}^c(\R^n)\to \mathrm{Spin}(\R^n) \to 
\mathrm{O}(n, \R)\) the class \(t_{\R^n}\) is the `Bott' class figuring in 
equivariant Bott Periodicity. 
 
\end{example}

Certain further
normalizations can be made in order to describe the groups 
\(\RK_{\Grd, X}^*(Y)\). A standard one is to replace the \(\Z/2\)-grading on 
\(\mathcal{E}\) by the standard even grading, so that \(\mathcal{E}\) consists 
of the sum of two copies of the same Hilbert module. This means that 
\(F\) can be taken to be of the form 
\(\begin{pmatrix} 0 & F_1^*\\F_1&0\end{pmatrix}\) and the conditions involving 
\(F\) are replaced by ones involving \(F_1\) and \(F_1^*\); we may as well 
replace \(F\) by \(F_1\). With this convention, the Fredholm conditions are that
 \(f (FF^*-1)\) and \(f(F^*F-1\) are compact for all \(f\in \CONT_0(X)\).  
 In other words, \(y \mapsto F_y\) takes essentially unitary values in \(\mathbb{B} (\Hils_y)\) 
 for all \(y\in Y\) and  the compact-operator-valued functions 
 \(FF^*-1\) and \(F^*F-1\) vanish at infinity along the fibres of \(\anchor\colon Y \to X\). 

The equivariant stabilization theorem for Hilbert modules 
 implies that we may take \(\Hils\) to have 
 the special form \(L^2(\Grd)^\infty\otimes_{\CONT_0(\Base)} \CONT_0(Y)\), 
 where \(L^2(\Grd)\) is the \(\Grd\)-equivariant right Hilbert \(\CONT_0(\Base)\)-
 Hilbert module defined using the Haar system of \(\Grd\), and the superscript 
 indicates the sum of countably many copies of \(L^2(\Grd)\). The corresponding 
 field of Hilbert space has value \(L^2(\Grd^y)^\infty\) at \(y\in Y\) where 
 \(\Grd^y\) denotes all \(G \in \Grd\) ending in \(y\), on which we have a 
 given measure specified by the Haar system of \(\Grd\).  
 
 This leads to a 
 description of \(\RK^0_{\Grd, Y} (X)\) as the 
group of homotopy-classes of \(\Grd\)-equivariant continuous maps from 
\(Y\) to the space \(\mathcal{F}_\Grd\) of Fredholm operators on the Hilbert 
spaces \(L^2(\Grd^y)^\infty\), but topologizing the space \(\mathcal{F}_\Grd\) 
is somewhat delicate. Similarly, the relative groups \(\RK^*_{\Grd,X}(Y)\) 
are maps to Fredholm operators with compact vertical support with 
respect to the map \(Y \to X\), where the 
\emph{support} of a map to Fredholm operators is by definition the complement of 
the set where the map takes invertible values. 

\begin{remark}
If \(\Grd\) acts properly and co-compactly on \(X\), \(A\) is a trivial 
\(\Grd\)-C*-algebra and \(B\) is a 
\(\Grd\ltimes X\)-C*-algebra, then 
there is a canonical isomorphism 
\[ \KK^{\Grd\ltimes X} (\CONT_0(X)\otimes A, B) \cong \KK(A, \Grd\ltimes B).\]
In particular, the \(\Grd\)-equivariant representable \(\K\)-theory of \(X\) 
agrees with the \(\K\)-theory of the corresponding cross-product. 
Under this identification, classes in \(\RK^0_\Grd (X)\) which are represented by 
equivariant vector bundles on \(X\) correspond to classes in 
\(\K_0(\Grd\ltimes \CONT_0(X))\) which are represented by 
projections in the stabilisation of \(\Grd\ltimes \CONT_0(X)\). 
See \cite{Emerson-Meyer:Equivariant_K} for more information.

\end{remark}

Thus, even if the reader is only interested in groups, or the trivial group, 
it is convenient to introduce groupoids to some extent in order to 
describe cohomology theories with different support conditions. 

\subsection{Tensor and forgetful functors }
The following simple functor will play an important role. 
If \(\Dual\) is a \(\Grd\ltimes \Tot\)-algebra, we denote by \(T_\Dual\) the map 
\[ \RKK^\Grd(\Tot; A,B) \defeq \KK^{\Grd\ltimes X} (\CONT_0(X)\otimes A, 
\UNIT_X\otimes B) \to \KK^\Grd(\Dual\otimes A, \Dual\otimes B)\]
which sends a \(\Grd\ltimes \Tot\)-equivariant right Hilbert 
\(\CONT_0(X)\otimes  B\)-Hilbert module \(\mathcal{E}\) to 
\(\mathcal{E}\otimes_X \Dual\), the tensor product being \emph{in the 
category of \(\Grd\ltimes X\)-algebras} (we accordingly use a subscript for 
emphasis) and sends 
\(F \in \mathbb{B}(\mathcal{E})\)
 to the operator \(F\otimes_X \ID_\Dual\). This definition 
 makes sense since \(F\) commutes with 
 the \(\CONT_0(\Tot)\)-structure on \(\mathcal{E}\). 
 
 The functor \(T_\Dual\) 
 is the composition of external product
 \[\blank  \otimes_X 1_\Dual \colon \KK^{\Grd\ltimes X}
 (\CONT_0(X)\otimes A, \CONT_0(X)\otimes B) 
 \to \KK^{\Grd\ltimes X} (A\otimes_X \Dual, B\otimes_X \Dual )\]
 (where the \(X\)-structure on \(A\otimes \Dual\) \emph{etc}. 
 is on the \(\Dual\) factor), and the forgetful map 
 \[ \KK^{\Grd\ltimes X} (A\otimes_X \Dual, B_X\otimes \Dual) \to 
 \KK^\Grd(A\otimes_X \Dual, B\otimes_X \Dual)\] 
which maps a \(\Grd\ltimes X\)-algebra or Hilbert module to 
the underlying \(\Grd\)-algebra, or Hilbert module, thus forgetting 
the \(X\)-structure.

\subsection{Kasparov duals}
\label{sec:abstract_duals}
We begin our discussion of duality by 
 by formalizing some duality calculations of Kasparov, \emph{c.f.} 
\cite{Kasparov:Novikov}*{Theorem 4.9}. Explicit examples will be discussed 
later. 

For convenience of notation we will often write 
\(\UNIT\defeq \CONT_0(\Base)\). This notation expresses the fact that
 \(\CONT_0(\Base)\) is the tensor
 unit in the tensor category of \(\Grd\)-C*-algebras. Similarly, if 
 \(\Grd\) acts on a space 
\(\Tot\) then we sometimes denote by  \(\UNIT_\Tot\) the \(\Grd\)-C*-algebra 
\(\CONT_0(\Tot)\); thus \(\UNIT_X\) is the tensor unit in the category 
of \(\Grd\ltimes X\) C*-algebras, \(\Tot\) being 
the base of 
\(\Grd\ltimes \Tot\). This notation is consistent with the source of this 
material (see \cite{Emerson-Meyer:Dualities}.)

\begin{definition}
  \label{def:Kasparov_dual}
  Let \(n\in\Z\).  An \emph{\(n\)\nb-dimensional \(\Grd\)\nb-equivariant Kasparov dual} for the \(\Grd\)\nb-space~\(\Tot\) is a triple \((\Dual,D,\Theta)\), where
  \begin{itemize}
  \item \(\Dual\) is a (possibly \(\Z/2\)-graded) \(\Grd\ltimes\Tot\)-\(\Cst\)\nb-algebra,

  \item \(D\in\KK^\Grd_{-n}(\Dual,\UNIT)\), and

  \item \(\Theta\in\RKK^\Grd_n(\Tot;\UNIT,\Dual)\),
  \end{itemize}
  subject to the following conditions:
  \begin{enumerate}[label=\textup{(\arabic{*})}]
  \item \(\Theta \otimes_\Dual D = \ID_\UNIT\) in \(\RKK^\Grd_0(\Tot;\UNIT,\UNIT)\);

  \item \(\Theta \otimes f = \Theta \otimes_\Dual T_\Dual(f)\) in \(\RKK^\Grd_{*+n}(\Tot;A,B\otimes\Dual)\) for all \(\Grd\)\nb-\(\Cst\)\nb-algebras \(A\) and~\(B\) and all \(f\in\RKK^\Grd_*(\Tot;A,B)\);

  \item \(T_\Dual(\Theta)\otimes_{\Dual\otimes\Dual} \Phi_\Dual = (-1)^n T_\Dual(\Theta)\) in \(\KK^\Grd_n(\Dual,\Dual\otimes\Dual)\), where~\(\Phi_\Dual\) is the flip automorphism on 
  \(\Dual\otimes\Dual\).
  \end{enumerate}
\end{definition}


The following theorem is proved in \cite{Emerson-Meyer:Dualities}. 

\begin{theorem}
  \label{thm:first_duality}
  Let \(n\in\Z\), let~\(\Dual\) be a \(\Grd\ltimes\Tot\)-\(\Cst\)\nb-algebra, \(D\in\KK^\Grd_{-n}(\Dual,\UNIT)\), and \(\Theta\in\RKK^\Grd_n(\Tot;\UNIT,\Dual)\).  Define two natural transformations
  \begin{alignat*}{2}
    \PD&\colon \KK^\Grd_{i-n}(\Dual\otimes A,B) \to \RKK^\Grd_i(\Tot;A,B), &\qquad
    f&\mapsto \Theta \otimes_\Dual f,\\
    \PD^*&\colon \RKK^\Grd_i(\Tot;A,B) \to \KK^\Grd_{i-n}(\Dual\otimes A,B), &\qquad g&\mapsto (-1)^{in} T_\Dual(g) \otimes_\Dual D,
  \end{alignat*}
  These two are inverse to each other if and only if \((\Dual,D,\Theta)\) is an \(n\)\nb-dimensional \(\Grd\)\nb-equivariant Kasparov dual for~\(\Tot\).  
  
\end{theorem}

\subsection{Abstract duals}
The reader may have noticed that the only place the \(\CONT_0(X)\)-structure 
on \(\Dual\) comes into play in the conditions listed in 
Definition \ref{def:Kasparov_dual}, and in the statement of  
Theorem \ref{thm:first_duality}, is via the functor \(T_\Dual\). In particular, if one has 
a Kasparov dual \(( \Dual, D, \Theta)\) and if one changes the 
\(\CONT_0(\Tot)\)-structure on 
\(\Dual\), for example by composing it with a \(\Grd\)-equivariant 
homeomorphism of \(\Tot\), then the 
map \(\PD\) of Theorem \ref{thm:first_duality} does not change; 
since by the theorem \(\PD^*\) is its inverse map, it 
would not change either, strangely, since its definition uses \(T_\Dual\). 
In fact it turns out that the functor \(T_\Dual\) can be reconstructed 
from \(\PD\) if one knows that \(\PD\) is an isomorphism. 
This is an important idea in connection 
with the Lefschetz map and suggests the following useful definition. 

\begin{definition}
  \label{def:abstract_dual}
  An \(n\)\nb-dimensional \emph{abstract dual} for~\(\Tot\) is a pair \((\Dual,\Theta)\), where~\(\Dual\) is a \(\Grd\)\nb-\(\Cst\)\nb-algebra and \(\Theta\in\RKK^\Grd_n(\Tot;\UNIT,\Dual)\), such that the map \(\PD\) defined as in Theorem~\ref{thm:first_duality} is an isomorphism for all \(\Grd\)-\(\Cst\)\nb-algebras \(A\) and~\(B\).
\end{definition}

This definition is shorter, and, as mentioned, is useful for theoretical reasons, 
but it seems like it should be difficult to check in 
practise. 

In any case, it is clear from 
Theorem~\ref{thm:first_duality} that a pair \((\Dual,\Theta)\) is an abstract 
dual if it is part of a Kasparov dual \((\Dual,D,\Theta)\).

\begin{proposition}
  \label{pro:dual_unique}
  An abstract dual for a space~\(\Tot\) is unique up to a canonical \(\KK^\Grd\)-equivalence if it exists, and even covariantly functorial in the following sense.

  Let \(\Tot\) and~\(\Other\) be two \(\Grd\)\nb-spaces and let \(f\colon \Tot\to \Other\) be a \(\Grd\)\nb-equivariant continuous map.  Let \((\Dual_\Tot,\Theta_\Tot)\) and \((\Dual_\Other,\Theta_\Other)\) be abstract duals for \(\Tot\) and~\(\Other\) of dimensions \(n_\Tot\) and~\(n_\Other\), respectively.  Then there is a unique \(\Dual_f\in\KK^\Grd_{n_\Other-n_\Tot}(\Dual_\Tot,\Dual_\Other)\) with \(\Theta_\Tot \otimes_{\Dual_\Tot} \Dual_f = f^*(\Theta_\Other)\).  Given two composable maps between three spaces with duals, we have \(\Dual_{f\circ g} = \Dual_f \circ \Dual_g\).  If \(\Tot=\Other\), \(f=\ID_\Tot\), and \((\Dual_\Tot,\Theta_\Tot)=(\Dual_\Other,\Theta_\Other)\), then \(\Dual_f=\ID_{\Dual_\Tot}\).  If only \(\Tot=\Other\), \(f=\ID_\Tot\), then~\(\Dual_f\) is a \(\KK^\Grd\)-equivalence between the two duals of~\(\Tot\).
\end{proposition}

Although the map \(f\colon X \to Y\) appearing in Proposition \ref{pro:dual_unique} does not have 
to be proper, it nonetheless yields a morphism \(\Dual_f\) 
in \(\KK^\Grd\).

\subsection{Duality co-algebra}
\label{subsec:duality_coalgebra}
Let \((\Dual,\Theta)\) be an \(n\)\nb-dimensional abstract dual for a \(\Grd\)\nb-space~\(\Tot\). 
 By the Yoneda Lemma, another abstract dual \((\Dual',\Theta')\) also for \(X\) and say 
of dimension~\(n'\) is related to \((\Dual,\Theta)\) by an invertible
 element
 \begin{equation}
 \label{eq:change_of_dual}
 \psi\in\KK^\Grd_{n'-n}(\Dual,\Dual'), \;\textup{such that}\;
 \Theta\otimes_{\Dual} \psi = \Theta'.\end{equation}  
We repeat for emphasis that since \( (\Dual, \Theta)\) is
 only an \emph{abstract} dual, 
we are not assuming 
that there is a \(\Grd\ltimes X\)-structure on \(\Dual\).
However, we are going to attempt to reconstruct what we might consider to be 
a \(\Grd\ltimes X\)-structure on \(\Dual\) \emph{at the level of \(\KK\)-theory}.
 Along the way we will keep track of how 
the change in dual from \((\Dual, \Theta)\) to \((\Dual', \Theta')\) 
affects our constructions. 

Define \(D\in\KK^\Grd_{-n}(\Dual,\UNIT)\) by the requirement 
\begin{equation}
  \label{eq:D_and_Theta}
  \PD(D)\defeq \Theta\otimes_\Dual D = 1_\UNIT
  \qquad \text{in \(\RKK^\Grd_0(\Tot;\UNIT,\UNIT)\).}
\end{equation}
as in the first condition in Definition~\ref{def:Kasparov_dual}). 
The class \(D\) should thus play the role of 
 the class named~\(D\) in a Kasparov dual. 
 It is routine to check that when 
 we change the dual, as above, \(D\) is replaced by \(\psi^{-1}\otimes_\Dual D\).

We call~\(D\) \emph{counit of the duality} because it plays the algebraic role of a counit in the 
theory of adjoint functors (see \cite{Emerson-Meyer:Dualities} and also Remark 
\ref{subsec:intro_to_embedding_theorems} below).

Define \(\comul\in\KK^\Grd_n(\Dual,\Dual\otimes\Dual)\) by the requirement that 
\[
\PD(\comul) \defeq \Theta\otimes_\Dual \comul = \Theta \otimes_\Tot \Theta \qquad \text{in \(\RKK^\Grd_{2n}(\Tot;\UNIT,\Dual\otimes\Dual)\).}
\]
We call~\(\comul\) the \emph{comultiplication of the duality}.  When we change the dual, \(\comul\) is replaced by
\[
(-1)^{n(n'-n)}\psi^{-1}\otimes_\Dual \comul \otimes_{\Dual\otimes\Dual} (\psi\otimes\psi) \in \KK^\Grd_{n'}(\Dual',\Dual'\otimes\Dual').
\]

\begin{remark}
\label{rem:duality_coalgebra}
If \(n=0\) then the object~\(\Dual\) of \(\KK^\Grd\) with counit~\(D\) and comultiplication~\(\comul\) is a cocommutative, counital coalgebra object in the tensor category \(\KK^\Grd\): 
  \begin{gather}
    \label{eq:coassociative}
     \comul\otimes_{\Dual\otimes \Dual}(\comul\otimes 1_\Dual) = \comul \otimes_{\Dual\otimes\Dual}(1_\Dual\otimes\comul),
    \\
    \label{eq:cocommutative}
    \comul\otimes_{\Dual\otimes\Dual}\Phi_\Dual = \comul,
    \\
    \label{eq:counit}
    \comul \otimes_{\Dual\otimes\Dual} (D\otimes 1_\Dual) = 1_\Dual = \comul \otimes_{\Dual\otimes\Dual} (1_\Dual\otimes D).
  \end{gather}
  Equation~\eqref{eq:coassociative} holds in \(\KK^\Grd_{2n}(\Dual,\Dual^{\otimes 3})\), equation~\eqref{eq:cocommutative} holds in \(\KK^\Grd_n(\Dual,\Dual\otimes\Dual)\), and~\eqref{eq:counit} holds in \(\KK^\Grd_0(\Dual,\Dual)\).
\end{remark}

Now, for \(\Grd\)\nb-\(\Cst\)\nb-algebras \(A\) and~\(B\), we define
\[
T'_\Dual\colon \RKK^\Grd_*(\Tot;A,B) \to \KK_*^\Grd(\Dual\otimes A,\Dual\otimes B), \qquad f\mapsto \comul \otimes_\Dual \PD^{-1}(f),
\]
where~\(\PD\) is the duality isomorphism, \(\comul\) is the comultiplication of the duality, and~\(\otimes_\Dual\) operates on the \emph{second} copy of~\(\Dual\) in the target \(\Dual\otimes\Dual\) of~\(\comul\).  A computation yields that 
\begin{equation}
  \label{eq:PD_sigma_prime}
  \PD\bigl(T'_\Dual(f)\bigr) = \Theta \otimes_\Tot f
  \qquad \text{in \(\RKK^\Grd_{i+n}(\Tot;A,\Dual\otimes B)\)}
\end{equation}
for all \(f\in\RKK^\Grd_i(\Tot;A,B)\).  It follows that 
\[
T'_\Dual(f) = T_\Dual(f)
\]
if \((\Dual,\Theta)\) is part of a Kasparov dual, and thus 
\(T_\Dual\) is in fact independent of the \(\Grd\ltimes X\)-structure on \(\Dual\), 
verifying our guess above. 

When we change the dual, we replace~\(T'_\Dual\) by the map
\begin{equation}
\label{eq:functoriality_of_T}
\RKK^\Grd_i(\Tot;A,B) \ni f\mapsto (-1)^{i(n-n')}\psi^{-1}\otimes_\Dual T_\Dual(f)\otimes_\Dual \psi \in \KK^\Grd_i(\Dual'\otimes A, \Dual'\otimes B).
\end{equation}

In fact, one can check that the maps \(T'_\Dual\) above define a functor
  \[
  T'_\Dual\colon \RKK^\Grd(\Tot)\to\KK^\Grd.
  \]
  This is a \(\KK^\Grd\)-functor in the sense that 
  it is compatible with the tensor products~\(\otimes\), and it is left adjoint to the functor \(p_\Tot^*\colon \KK^\Grd\to\RKK^\Grd\) induced from the groupoid homomorphism 
  \(\Grd\ltimes \Tot \to \Grd\). 

It follows that we can write the inverse duality 
map involved in an abstract dual \((\Dual, \Theta)\) as:
\begin{equation}
  \label{eq:inverse_PD_abstract}
  \PD^{-1}(f) = (-1)^{in} T'_\Dual(f) \otimes_\Dual D
  \qquad
  \text{in \(\KK^\Grd_{i-n}(\Dual\otimes A,B)\)}
\end{equation}
for \(f\in \RKK^\Grd_i(\Tot;A,B)\).
By the above discussion this formula agrees with the 
map \(\PD^*\) when we have a Kasparov dual. 

\subsection{The Lefschetz map }
\label{subsec:the_lefschetz_map}
The formal computations summarized in the previous section allows us to 
single out an interesting invariant of a 
\(\Grd\)-space \(X\), at least under the hypothesis that \(X\) has \emph{some}
abstract dual. 

For any \(\Grd\)-space \(X\) the
 diagonal embedding \(\Tot\to\Tot\times_\Base\Tot\) is a proper map and hence induces a \(^*\)\nb-homomorphism
\[
\UNIT_\Tot\otimes \UNIT_\Tot \cong \CONT_0(\Tot\times_\Base\Tot)\to\CONT_0(\Tot) = \UNIT_\Tot.
\]
This map is \(\Grd\ltimes\Tot\)-equivariant and hence yields
\[
\Delta_\Tot\in\RKK^\Grd(\Tot;\UNIT_\Tot,\UNIT) \cong \KK^{\Grd\ltimes\Tot}\bigl(\CONT_0(\Tot\times_\Base\Tot),\CONT_0(\Tot)\bigr).
\]
We call this the \emph{diagonal restriction class}. It yields a canonical map
\begin{equation}
  \label{eq:diagonal_gives_this}
  \blank\otimes_{\UNIT_\Tot}\Delta_\Tot\colon
  \KK^\Grd(\UNIT_\Tot\otimes A,\UNIT_\Tot\otimes B) \to
  \RKK^\Grd(\Tot;\UNIT_\Tot\otimes A,B).
\end{equation}
In particular, this contains a map \(\KK^\Grd(\UNIT_\Tot,\UNIT_\Tot) \to \RKK^\Grd(\Tot;\UNIT_\Tot,\UNIT)\).

\begin{example}
  \label{exa:diagonal_restrict_proper_map}
  If \(f\colon \Tot\to\Tot\) is a proper, continuous, \(\Grd\)\nb-equivariant map, then
  \[
  [f]\otimes_{\UNIT_\Tot} \Delta_\Tot \in \RKK^\Grd(\Tot;\UNIT_\Tot,\UNIT)
  \]
  is the class of the \(^*\)\nb-homomorphism induced by \((\ID_\Tot,f)\colon \Tot \to \Tot \times_\Base\Tot\).

  Now drop the assumption that~\(f\) be proper.  Then \((\ID_\Tot,f)\) is still a proper, continuous, \(\Grd\)\nb-equivariant map.  The class of the \(^*\)\nb-homomorphism it induces is equal to \(f^*(\Delta_\Tot)\), where we use the maps
  \[
  f^*\colon \RKK^\Grd_*(\Tot;A,B) \to \RKK^\Grd_*(\Tot;A,B)
  \]
  for \(A=\UNIT_\Tot\), \(B=\UNIT\) induced by \(f\colon \Tot\to\Tot\) (the functor 
  \( X\mapsto \RKK^\Grd (X ; A,B)\) is 
  functorial with respect to arbitrary \(\Grd\)-maps, not just proper ones.) This suggests that we can think of \(\RKK^\Grd(X; \UNIT_X, \UNIT)\) as generalized, possibly \emph{non-proper} self-maps of 
  \(X\). 
\end{example}

In fact if the anchor map \(X\to \Base\) is a proper map, so that \(X\) is a bundle of 
compact spaces over \(\Base\), then  
\(\blank\otimes_{\UNIT_\Tot}\Delta_\Tot\) is an isomophism (an easy exercise in the 
definitions.) 

Now let~\(T'_\Dual\) be the tensor functor and~\(\Delta_\Tot\) the diagonal restriction class of an abstract dual.  We define the \emph{multiplication class} of~\(\Dual\) by
\begin{equation}
  \label{eq:def_multiplication_class}
  [m] \defeq T'_\Dual(\Delta_\Tot)
  \in \KK^\Grd_0(\Dual\otimes\UNIT_\Tot,\Dual).
\end{equation}
A change of dual as in \eqref{eq:change_of_dual} replaces \([m]\) by \(\psi^{-1}\otimes_\Dual [m]\otimes_\Dual \psi\).

\begin{lemma}
  \label{lem:multiplication_class}
  Let \((\Dual,D,\Theta)\) be a Kasparov dual.  Then~\([m]\) is the class in \(\KK^\Grd\) of the multiplication homomorphism \(\CONT_0(\Tot)\otimes_\Base \Dual\to\Dual\) that describes the \(\Tot\)\nb-structure on~\(\Dual\) \textup{(}up to commuting the tensor factors\textup{)}.
\end{lemma}

We now have enough theoretical development to define the Lefschetz map and 
sketch the proof of its homotopy invariance.

Let~\(\Tot\) be a \(\Grd\)\nb-space and \((\Dual,\Theta)\) an \(n\)\nb-dimensional abstract dual for~\(\Tot\), \(\PD\) and \(\PD^{-1}\) the duality isomorphisms.  As before, we write
\(
\UNIT \defeq \CONT_0(\Base)\), \(\UNIT_\Tot \defeq \CONT_0(\Tot)
\)
and  \(\Delta_\Tot\in \RKK^\Grd(\Tot;\UNIT_\Tot,\UNIT) = \KK^{\Grd\ltimes\Tot}(\UNIT_\Tot\otimes\UNIT_\Tot,\UNIT_\Tot)\) the diagonal restriction class and
\[
\bar{\Theta} \defeq \forget_\Tot(\Theta) \in \KK^\Grd_n(\UNIT_\Tot,\Dual\otimes\UNIT_\Tot).
\]

\begin{definition}
  \label{def:Lefschetz_map}
  The equivariant \emph{Lefschetz map}
  \[
  \Lef\colon \RKK^\Grd_*(\Tot;\UNIT_\Tot,\UNIT) \to \KK_*^\Grd(\UNIT_\Tot,\UNIT)
  \]
  for a \(\Grd\)\nb-space~\(\Tot\) is defined as the composite map
  \[
  \RKK^\Grd_i(\Tot;\UNIT_\Tot,\UNIT) \xrightarrow{\PD^{-1}} \KK^\Grd_{i-n}(\Dual\otimes\UNIT_\Tot,\UNIT) \xrightarrow{\bar{\Theta}\otimes_{\Dual\otimes\UNIT_\Tot}\blank} \KK^\Grd_i(\UNIT_\Tot,\UNIT).
  \]
  The equivariant \emph{Euler characteristic} of~\(\Tot\) is
  \[
  \Eul_\Tot \defeq \Lef(\Delta_\Tot) \in \KK^\Grd_0(\UNIT_\Tot,\UNIT) = \KK^\Grd_0\bigl(\CONT_0(\Tot),\CONT_0(\Base)\bigr).
  \]
\end{definition}

Let \(f\in \RKK^\Grd_i(\Tot;\UNIT_\Tot,\UNIT)\).  Equations \eqref{eq:inverse_PD_abstract} and~\eqref{eq:def_multiplication_class} yield
\begin{align}
  \label{eq:compute_Lef}
  \Lef(f) &= (-1)^{in} \bar{\Theta}
  \otimes_{\Dual\otimes\UNIT_\Tot} T'_\Dual(f) \otimes_\Dual D,\\
  \label{eq:compute_Eul}
  \Eul_\Tot &= (-1)^{in} \bar{\Theta}
  \otimes_{\Dual\otimes\UNIT_\Tot} [m] \otimes_\Dual D.
\end{align}

We have already established that 
if \((\Dual,\Theta)\) is part of a Kasparov dual, then \(T'_\Dual=T_\Dual\) and \([m]\) is the \(\KK\)-class of the multiplication \(^*\)\nb-homomorphism \(\CONT_0(X,\Dual)\to\Dual\), so that~\eqref{eq:compute_Lef} yields explicit formulas for \(\Lef(f)\) and \(\Eul_X\).  This is extremely 
important because otherwise it would not be possible to compute these invariants. 

Let \(\Tot\) and~\(\Tot'\) be \(\Grd\)\nb-spaces, and let \(f\colon \Tot\to\Tot'\) be a \(\Grd\)\nb-homotopy equivalence.  Then~\(f\) induces an equivalence of categories \(\RKK^\Grd(\Tot') \cong \RKK^\Grd(\Tot)\), that is, we get invertible maps
\[
f^*\colon \RKK^\Grd_*(\Tot';A,B) \to \RKK^\Grd_*(\Tot;A,B)
\]
for all \(\Grd\)\nb-\(\Cst\)\nb-algebras \(A\) and~\(B\).  Now assume, in addition, that~\(f\) is proper; we do not need the inverse map or the homotopies to be proper.  Then~\(f\) induces a \(^*\)\nb-homomorphism \(f^!\colon \CONT_0(\Tot') \to \CONT_0(\Tot)\), which yields \([f^!]\in \KK^\Grd\bigl(\CONT_0(\Tot'),\CONT_0(\Tot)\bigr)\).  We write \([f^!]\) instead of \([f^*]\) to better distinguish this from the map~\(f^*\) above.  Unless~\(f\) is a \emph{proper} \(\Grd\)\nb-homotopy equivalence, \([f^!]\) need not be invertible.

\begin{theorem}
  \label{thm:Lef_well-defined}
  Let \(\Tot\) and~\(\Tot'\) be \(\Grd\)\nb-spaces with abstract duals, and let \(f\colon \Tot\to\Tot'\) be both a proper map and a \(\Grd\)\nb-homotopy equivalence.  Then
  \[
  [f^!]\otimes_{\CONT_0(\Tot)} \Eul_\Tot = \Eul_{\Tot'} \qquad \text{in \(\KK^\Grd_0(\CONT_0(\Tot'),\UNIT)\)}
  \]
  and the Lefschetz maps for \(\Tot\) and~\(\Tot'\) are related by a commuting diagram
  \[
  \xymatrix{ \RKK^\Grd_*(\Tot;\CONT_0(\Tot),\UNIT) \ar[d]^{\Lef_\Tot} & \ar[l]_{f^*}^{\cong} \RKK^\Grd_*(\Tot';\CONT_0(\Tot),\UNIT) \ar[r]^{[f^!]^*} &
    \RKK^\Grd_*(\Tot';\CONT_0(\Tot'),\UNIT) \ar[d]^{\Lef_{\Tot'}} \\
    \KK_*^\Grd(\CONT_0(\Tot),\UNIT) \ar[rr]^{[f^!]^*} && \KK_*^\Grd(\CONT_0(\Tot'),\UNIT), }
  \]
  where \([f^!]^*\) denotes composition with~\([f^!]\).

  In particular, \(\Eul_\Tot\) and the map \(\Lef_\Tot\) do not depend on the chosen dual.
\end{theorem}

The proof relies on the discussion preceding the theorem. 

Theorem~\ref{thm:Lef_well-defined} implies that the Lefschetz maps for properly \(\Grd\)\nb-homotopy equivalent spaces are equivalent because then \([f^!]\) is invertible, so that all horizontal maps in the diagram in Theorem~\ref{thm:Lef_well-defined} are invertible.  In this sense, the Lefschetz map and the Euler class are invariants of the proper \(\Grd\)\nb-homotopy type of~\(\Tot\).

The construction in Example~\ref{exa:diagonal_restrict_proper_map} associates a class \([\Delta_f] \in \RKK^\Grd_0(\Tot;\CONT_0(\Tot),\UNIT)\) to any continuous, \(\Grd\)\nb-equivariant map \(f\colon \Tot\to\Tot\); it does not matter whether~\(f\) is proper.  We abbreviate
\[
\Lef(f) \defeq \Lef([\Delta_f])
\]
and call this the Lefschetz invariant of~\(f\).  Of course, equivariantly homotopic self-maps induce the same class in \(\RKK^\Grd(\Tot;\CONT_0(\Tot),\UNIT)\) and therefore have the same Lefschetz invariant.  We have \(\Lef(\ID_\Tot)=\Eul_\Tot\).

More generally, specializing \eqref{eq:diagonal_gives_this} gives a map 
\[
\blank\otimes_{\UNIT_\Tot}\Delta_\Tot\colon \KK_*^\Grd\bigl(\CONT_0(\Tot),\CONT_0(\Tot)\bigr) \to \RKK^\Grd_*(\Tot;\CONT_0(\Tot),\UNIT),
\]
which we compose with the Lefschetz map; abusing notation, we still denote this 
composition by 
\[
\Lef\colon \KK_*^\Grd\bigl(\CONT_0(\Tot),\CONT_0(\Tot)\bigr) \to \KK_*^\Grd(\CONT_0(\Tot),\UNIT)
\]

Finally, we record that 
Lefschetz invariants for elements of \(\RKK^\Grd_*(\Tot;\CONT_0(\Tot),\UNIT)\) can be arbitrarily complicated: the Lefschetz map is rather easily seen to be split surjective. 
The splitting is given by specializing the \emph{inflation map}
\begin{equation}
\label{eq:inflation_map}
p_X^*\colon \KK_*^\Grd(A, B) \to \KK^{\Grd\ltimes X}(\UNIT_X\otimes A, \UNIT_X\otimes B)
\end{equation}
to \(A \defeq \UNIT_X\) and \(B \defeq \UNIT\). 

The fundamental example of a Kasparov dual is provided by the vertical tangent 
space to a bundle of smooth manifolds over the base \(\Base\) of a groupoid, 
in which morphisms act smoothly. We come back to this in \S 
\ref{sec:dualities_and_lef}.

\section{Examples of computations of the Lefschetz map}
\label{subsec:map_examples}

In this section we will give some examples of computations of the 
Lefschetz map for various instances of spaces with 
duals and for equivariant self-morphisms coming 
from actual \emph{maps}. The problem of computing the 
Lefschetz invariants of more general Kasparov self-morphisms for 
the next section is a central problem for us, and will be treated 
in \S \ref{sec:dualities_and_lef} once we have available the
theory of equivariant correspondences.  

Most of the examples are quite close to 
proper actions, but they do not quite have to be proper. The point is that an abstract 
(or Kasparov) dual for a \(\Grd\)-space \(X\) also yields one for \(X\) regarded as a 
\(\Grd'\)-space where \(\Grd'\) is a (not necessarily closed) subgroupoid of 
\(\Grd\). Even if \(X\) is proper as a \(\Grd\)-space, it need not be proper as a 
\(\Grd'\)-space since \(\Grd'\) may not be closed. 

Known examples of duals are of 
 two types; they are both Kasparov duals. 
 If \(\Grd\) is a locally compact group acting as simplicial automorphisms of a 
finite-dimensional 
simplicial complex then \(X\) has a Kasparov dual by \cite{Emerson-Meyer:Euler}. 
This does not quite imply that the action of \(\Grd\) is proper, since the action of 
\(\Grd\) could be trivial. Neither does it quite imply that \(\Grd\) must be discrete, 
though since the 
connected component of the identity of \(\Grd\) must act trivially, the only interesting 
examples must involve disconnected groups. For 
instance, the non-discrete group \(\mathrm{SL}(2, \Q_p)\) acts simplicially on a tree. 

 If \(X\) is a complete 
complete Riemannian manifold then \(X\) 
is a proper \(\Grd\)-space where \(\Grd\) is the Lie group of isometries of \(X\),
 and either the Clifford algebra of 
\(X\) or the C*-algebra of \(\CONT_0\)-functions on the tangent bundle \(\tang\) 
of \(X\) is part of a Kasparov dual for \(X\) (see \emph{e.g.}
\cite{Emerson-Meyer:Euler} or 
\cite{Emerson-Meyer:Dualities}), and also \S \ref{sec:dualities_and_lef}. 
Hence if \(\Grd\) is \emph{any} group
of isometries of a Riemannian manifold \(X\) with finitely many components, 
then the \(\Grd\)-space \(X\) also has a Kasparov dual. If the tangent bundle 
\(\tang\) admits an equivariant \(\K\)-orientation, then \(\CONT_0(\tang)\) can 
be replaced by \(\CONT_0(X)\). For instance, the circle with 
the group \(\Z\) acting by an irrational rotation has a Kasparov dual of of dimension 
\(1\), given 
by \( (\CONT(\mathbb{T}), D, \Theta)\) where \(D\) is the class of the 
Dirac operator on the circle.

The Lefschetz invariants of equivariant self-maps of \(X\) will in both these 
situations turn out to be in some sense 
\emph{zero-dimensional} in the sense that they are built out of point-evaluation classes. 
As we will see later, the Lefschetz invariants of more general Kasparov 
self-morphisms are more complicated, higher-dimensional objects.

\subsection{The combinatorial Lefschetz map}
\label{subsec:combinatorial_case}
Let~\(X\) be a finite-dimensional simplicial complex and
let~\(G\) be a locally compact group acting smoothly and
simplicially on~\(X\) (that is, stabilisers of points are
open).   We follow \cite{Emerson-Meyer:Euler} and 
\cite{Emerson-Meyer:Equi_Lefschetz}. 
Assume that~\(X\) admits a colouring (that
is, \(X\) is typed) and that~\(G\) preserves the colouring.
This ensures that if \(g\in G\) maps a simplex to itself, then
it fixes that simplex pointwise.

Let \(SX\) be the set of (non-degenerate) simplices of~\(X\)
and let \(S_dX\subseteq SX\) be the subset of
\(d\)\nb-dimensional simplices.  The group~\(G\) acts
on the discrete set~\(SX\) preserving the
decomposition \(SX=\bigsqcup S_d X\).  Decompose~\(SX\) into
\(G\)\nb-orbits.  For each orbit \(\dot\sigma \subseteq SX\),
choose a representative \(\sigma\in SX\) and let
\(\xi_\sigma\in X\) be its barycentre and
\(\Stab(\sigma)\subseteq G\) its stabiliser.  Restriction to
the orbit~\(\dot\sigma\) defines a \(G\)\nb-equivariant
\Star{}homomorphism
\begin{equation}
  \label{eq:dot_xi}
  \xi_{\dot\sigma}\colon
  \CONT_0(X)\to\CONT_0\bigl(G/\Stab(\sigma)\bigr) \to
  \Comp\bigl(\ell^2(G/\Stab\sigma)\bigr),
\end{equation}
where the second map is the representation by pointwise
multiplication operators.  We let \([\xi_{\dot\sigma}]\) be its
class in \(\KK^G_0(\CONT_0(X),\UNIT)\).

Let \(\varphi\colon X\to X\) be a \(G\)\nb-equivariant
self-map of~\(X\).  Since~\(\varphi\) is
\(G\)\nb-equivariantly homotopic to a \(G\)\nb-equivariant
cellular map, we may assume without loss of generality
that~\(\varphi\) is itself cellular.  Hence it induces a
\(G\)\nb-equivariant chain map
\[
\varphi \colon C_\bullet(X)\to C_\bullet(X),
\]
where \(C_\bullet (X)\) is the chain complex of \emph{oriented}
simplices of~\(X\).  A basis for \(C_\bullet (X)\) is given by
the set of (un-)oriented simplices, by arbitrarily choosing an
orientation on each simplex.  We may describe the chain
map~\(\varphi\) by its matrix coefficients
\(\Selfmap_{\sigma\tau}\in\Z\) with respect to this basis; thus
the subscripts are unoriented simplices.  For example,
if~\(\Selfmap\) maps a simplex to itself, and reverses
orientation, then \(\Selfmap_{\sigma,\sigma}=-1\).
Since~\(\Selfmap\) is \(G\)\nb-equivariant,
\(\Selfmap_{g(\sigma),g(\sigma)} = \varphi_{\sigma\sigma}\).
So the following makes sense.

\begin{notation}
  \label{note:mult_gamma}
  For \(\dot{\sigma} \in G\backslash S_dX\), let
  \(\mult(\Selfmap,\dot\sigma) \defeq
  (-1)^d\Selfmap_{\sigma\sigma} \in \Z\) for any choice of
  representative \(\sigma \in \dot{\sigma}\).
\end{notation}

The following theorem is proved in \cite{Emerson-Meyer:Equi_Lefschetz} using 
the simplicial dual 
developed in \cite{Emerson-Meyer:Euler}, and inspired by 
ideas of Kasparov and Skandalis in \cite{Kasparov-Skandalis:Buildings}.

\begin{theorem}
  \label{the:Lef_combinatorial}
  Let~\(X\) be a finite-dimensional coloured simplicial
  complex and let~\(G\) be a locally compact group that acts
  smoothly and simplicially on~\(X\), preserving the colouring.
  Let \(\selfmap\colon X\to X\) be a \(G\)\nb-equivariant
  self-map.  Define \(\mult(\Selfmap,\dot\sigma)\in\Z\) and
  \([\xi_{\dot\sigma}]\in \KK^G(\CONT_0(X),\UNIT)\) for
  \(\dot\sigma\in G\backslash SX\) as above.  Then
  \[
  \Lef(\selfmap) = \sum_{\dot\sigma\in G\backslash SX}
  \mult(\Selfmap,\dot\sigma)[\xi_{\dot\sigma}].
  \]
\end{theorem}

In the non-equivariant situation, if \(\Tot\) is connected and has only a 
finite number of simplices, then 
the formula just given reduces 
to the ordinary Lefschetz number of \(G\) by the (standard) 
argument that proves that the Euler characteristic of a finite simplicial complex 
can be computed either by counting ranks of simplicial homology groups, 
or by directly counting the number of simplices in the complex. In the 
noncompact case this doesn't make sense anymore since the number of 
orbits of simplices may be infinite, but our definition of the Lefschetz invariant
still makes sense. This illustrates the advantage of considering 
\(\K\)-homology classes instead of numbers as fixed-point data.

\subsection{The smooth Lefschetz map}
\label{sec:orient_vb_auto}
Let \(X\) be
 a smooth manifold and \(G\) a group acting by isometric 
 diffeomorphisms of \(X\). Then 
we can build a Kasparov dual for \(X\) using the Clifford algebra. The real (resp. 
complex) Clifford 
algebra of a Euclidean vector space \(V\) is generated as a real (resp. complex)
 unital *-algebra by 
an orthonormal basis \(\{e_i\}\) for \(V\) with the relations that \(e_i\) is self-adjoint 
and \(e_ie_j+e_je_i= 2\delta_{ij}\). In the complex case 
this produces a finite-dimensional \(\Z/2\)-graded
C*-algebra. More generally if \(V\) is a Euclidean vector bundle, this construction 
applies and produces a locally trivial bundle of finite-dimensional \(\Z/2\)-graded
C*-algebras over \(X\). If \(V\defeq \Tvert \Tot\) for a Riemannian manifold \(X\) then 
the \emph{Clifford algebra} of \(X\) is the corresponding C*-algebra of sections
vanishing at infinity. It is denoted \(\ctau (X)\). 
This C*-algebra carries a canonical action of the group of isometries of 
\(X\) and hence likewise for any subgroup. 

We discuss the specific mechanics of the Clifford Kasparov dual to the following 
extent. Let \(d\) be the de Rham differential on \(X\), acting on \(L^2\)-forms 
on \(X\). This Hilbert space carries a unitary action of \(G\) and both \(d\) and its 
adjoint \(d^*\) are \(G\)-equivariant, and \(d+d^*\) is an elliptic operator called the 
\emph{Euler} (or \emph{de Rham}) operator on \(X\).  If 
\(\omega\) is a differential form on \(X\) vanishing at infinity 
then the operator \(\lambda_\omega\) 
of exterior product with \(\omega\) defines an operator on \(L^2\)-forms which 
is bounded and the assignment 
\(\omega \mapsto \lambda_\omega + \lambda_\omega^*\) determines a 
representation of \(\ctau(X)\) which graded commutes modulo bounded 
operators with \(d+d^*\). Hence we get a cycle for \(\KK^G(\ctau(X), \C)\).
It represents the class \(D\) appearing in the Kasparov dual.

We start by recording one of the easiest computations of the Lefschetz map. 
The interested reader 
can easily prove it for his or herself using after looking briefly at the definition 
of the class \(\Theta\) (see \cite{Emerson-Meyer:Euler}) and reviewing the 
definition of the Lefschetz map. Recall that the Euler class of \(X\) is the 
Lefschetz invariant of the identity map of \(X\). 

\begin{proposition}
\label{pro:euler}
The Euler class of \(X\) is the class of the Euler operator on \(X\). 
\end{proposition}

The functoriality result \ref{thm:Lef_well-defined} combines with 
Proposition \ref{pro:euler} to imply the homotopy-invariance of the 
class of the Euler operator: 
\[ f_*(\Eul_X) = \Eul_{X'}\] for a proper 
\(\Grd\)-equivariant homotopy-equivalence \(f\colon X \to X'\) between 
smooth and proper \(\Grd\)-manifolds.



We now describe the Lefschetz invariant of a more general 
smooth \(G\)-equivariant self-map of \(X\). This requires a preliminary discussion. 

Let~\(Y\) be a locally
compact space and~\(G\) be a locally compact group
 acting continuously on~\(Y\), and let \(\pi\colon
E\to Y\) be a \(G\)\nb-equivariant Euclidean
\(\R\)\nb-vector bundle over~\(E\). 
Let \(A\colon E\to E\) be a \(G\)\nb-equivariant
vector bundle automorphism, that is, a continuous map \(E\to
E\) over~\(Y\) that restricts to \(\R\)\nb-vector space
isomorphisms on the fibres of~\(E\).  We are going to define a
\(G\)\nb-equivariant \(\Z/2\)-graded real line bundle
\(\sign(A)\) over~\(Y\). (Since we work with complex \(\K\)-theory
we will only use its complexification.) 

If~\(Y\) is a point, then \(G\)\nb-equivariant 
real vector bundles over~\(Y\) correspond to real orthogonal
representations of~\(G\). The endomorphism \(A\) 
becomes in this case an invertible linear map \(A\colon
\R^n\to\R^n\) commuting with~\(G\).  The sign is a virtual
\(1\)\nb-dimensional representation of~\(G\) and hence
corresponds to a pair \((\chi,n)\), where \(n\in\{0,1\}\) is
the parity (we are referring to the grading, either even or odd) 
of the line bundle and \(\chi\colon G\to \{-1,+1\}\)
is a real-valued character.  The overall parity will turn out
to be~\(0\) if~\(A\) preserves orientation and~\(1\) if~\(A\)
reverses orientation (see Example \ref{exa:orientation_representation}). 
 In this sense, our invariant will refine the orientation of~\(A\).  

  As above, let \(\Cliff(E)\) be the bundle of \emph{real} 
  Clifford algebras
  associated to~\(E\). We can also define in an analogous way 
  \(\Cliff(E)\) if~\(E\) carries an indefinite
bilinear form and it is a well-known fact from algebra 
that if the index of the bilinear form on~\(E\) is
divisible by~\(8\), then the fibres of \(\Cliff(E)\) are
isomorphic to matrix algebras.  In this case, a
\emph{\(G\)\nb-equivariant spinor bundle} for~\(E\) is a
\(\Z/2\)-graded real vector bundle~\(S_E\) together with a
grading preserving, \(G\)\nb-equivariant \Star{}algebra
isomorphism \(c\colon \Cliff(E)\to \operatorname{End}(S_E)\).
This representation is determined uniquely by its restriction
to \(E\subseteq \Cliff(E)\), which is a \(G\)\nb-equivariant
map \(c\colon E\to \operatorname{End}(S_E)\) such that \(c(x)\)
is odd and symmetric and satisfies \(c(x)^2 = \norm{x}^2\) for
all \(x\in E\).  

The spinor bundle is unique up to tensoring with a
\(G\)\nb-equivariant real line bundle~\(L\): if \(c_t\colon
E\to S_t\) for \(t=1,2\) are two \(G\)\nb-equivariant spinor
bundles for~\(E\), then we define a \(G\)\nb-equivariant real
line bundle~\(L\) over~\(Y\) by
\[
L\defeq \Hom_{\Cliff(E)}(S_1,S_2),
\]
and the evaluation isomorphism \(S_1\otimes L
\xrightarrow{\cong} S_2\) intertwines the representations
\(c_1\) and~\(c_2\) of \(\Cliff(E)\).

\begin{definition}
  \label{def:sign_A_K-oriented}
  Let \(A\colon E\to E\) be a real \(G\)\nb-equivariant vector
  bundle automorphism and let \(A=T\circ
  (A^*A)^{\nicefrac{1}{2}}\) be its polar decomposition with an
  orthogonal vector bundle automorphism \(T\colon E\to E\).

  Let~\(F\) be another \(G\)\nb-equivariant vector bundle
  over~\(Y\) with a non-degenerate bilinear form, such that the
  signature of \(E\oplus F\) is divisible by~\(8\), so that
  \(\Cliff(E\oplus F)\) is a bundle of matrix algebras
  over~\(\R\).  We assume that \(E\oplus F\) has a
  \(G\)\nb-equivariant spinor bundle, that is, there exists a
  \(G\)\nb-equivariant linear map \(c\colon E\oplus F\to
  \operatorname{End}(S)\) that induces an isomorphism of graded
  \Star{}algebras
  \[
  \Cliff(E\oplus F) \cong \operatorname{End}(S).
  \]
  Then
  \[
  c'\colon E\oplus F\to\operatorname{End}(S), \qquad
  (\xi,\eta)\mapsto c(T(\xi),\eta)
  \]
  yields another \(G\)\nb-equivariant spinor bundle for
  \(E\oplus F\).  We let
  \[
  \sign(A) \defeq
  \Hom_{\Cliff(E\oplus F)}\bigl((S,c'),(S,c)\bigr).
  \]
  This is a \(G\)\nb-equivariant \(\Z/2\)-graded real line
  bundle over~\(Y\).
\end{definition}
 
   It is not hard to check that \(\sign(A)\) is well-defined and  
   a homotopy invariant. Furthermore, \(\sign(A_1\circ A_2) \cong
    \sign(A_1)\otimes\sign(A_2)\) for two equivariant
    automorphisms \(A_1,A_2\colon E\rightrightarrows E\) of the
    same bundle, and \(\sign(A_1\oplus A_2) \cong
    \sign(A_1)\otimes\sign(A_2)\) for two equivariant vector
    bundle automorphisms \(A_1\colon E_1\to E_1\) and
    \(A_2\colon E_2\to E_2\).

     If \(Y\) is a point and \(\Grd\) is trivial, then \(\sign(A) = \R\) for 
     orientation-preserving~\(A\) and \(\sign (A) = \R^\op\) for orientation-reversing~\(A\), 
     as claimed above.

\begin{example}
  \label{exa:orientation_representation}
  Consider \(G=\Z/2\). 
  Let \(\tau\colon G\to\{1\}\) be the trivial character and let
  \(\chi\colon G\to \{+1,-1\}\) be the non-trivial character.
  Let \(\R_\chi\) denote the real representation of 
  \(G\) on \(\R\) with character \(\chi\). This can be considered
  trivially graded; let \(\R^\op_\chi\) denote the same 
  representation but with the opposite grading (the whole 
  vector space is considered odd.) 
   
   Consider \(A\colon \R_\chi\to\R_\chi\), \(t\mapsto -t\), so 
   \(A\) commutes with \(G\). 
  Then \(\sign(A)\cong\R_\chi^\op\) carries a non-trivial
  representation.

  To see this, let~\(F\) be~\(\R_\chi\) with
  negative definite metric.  Thus the Clifford algebra of
  \(\R_\chi\oplus\R_\chi\) is \(\Cliff_{1,1} \cong
  \mathbb{M}_{2\times2}(\R)\).  Explicitly, the map
  \[
  c(x,y)
  =\begin{pmatrix}0&x-y\\x+y&0\end{pmatrix}
  \]
  induces the isomorphism.  We equip~\(\R^2\) with the
  representation \(\tau\oplus\chi\), so that~\(c\) is
  equivariant.

  Twisting by~\(A\) yields another representation
  \[
  c'(x,y)\defeq c(-x,y) = Sc(x,y)S^{-1}
  \qquad\text{with} \quad S=S^{-1}=
  \begin{pmatrix}0&1\\-1&0\end{pmatrix}.
  \]
  Since~\(S\) reverses the grading and exchanges the
  representations \(\tau\) and~\(\chi\), it induces an
  isomorphism \((\R_\tau\oplus\R_\chi^\op)\otimes \R_\chi^\op
  \xrightarrow{\cong} \R_\tau\oplus\R_\chi^\op\).  Hence
  \(\sign(A) = \R_\chi^\op\).
\end{example}

Now let~\(X\) be a smooth Riemannian manifold and assume
that~\(G\) acts on~\(X\) isometrically and continuously.

Let \(\selfmap\colon X\to X\) be a \(G\)\nb-equivariant
self-map of~\(X\).  In order to write down an explicit local
formula for \(\Lef(\selfmap)\), we impose the following
restrictions on~\(\selfmap\):
\begin{itemize}
\item \(\selfmap\) is smooth;

\item the fixed point subset \(\Fix(\selfmap)\) of~\(\selfmap\)
  is a submanifold of~\(X\);

\item if \((p,\xi)\in TX\) is fixed by the derivative
  \(D\selfmap\), then~\(\xi\) is tangent to \(\Fix(\selfmap)\).
\end{itemize}
The last two conditions are automatic if~\(\selfmap\) is
isometric with respect to \emph{some} Riemannian metric (not
necessarily the given one) and this of course applies by averaging 
the given metric if \(\selfmap\) has finite
order.

In the simplest case, \(\selfmap\) and \(\ID_X\) are
transverse, that is, \(\ID-D\selfmap\) is invertible at each
fixed point of~\(\selfmap\); this implies that~\(\selfmap\) has
isolated fixed points.  

To describe the Lefschetz invariant, we abbreviate \(Y\defeq
\Fix(\selfmap)\).  This is a closed submanifold of~\(X\) by
assumption.  Let~\(\nu\) be the normal bundle of~\(Y\)
in~\(X\).  Since the tangent space of~\(Y\) is left fixed by
the derivative \(D\selfmap\), it induces a linear map
\(D_\nu\selfmap\colon \nu\to\nu\).  By assumption, the map
\(\ID_\nu-D_\nu\selfmap\colon \nu\to\nu\) is invertible.

\begin{theorem}
  \label{the:Lef_smooth}
  Let~\(X\) be a complete smooth Riemannian manifold, let~\(G\)
  be a locally compact group that acts on~\(X\) smoothly 
  and by isometries, and let \(\selfmap\colon X\to X\) be a
  self-map satisfying the three conditions enumerated above. 
  Let~\(\nu\) be the normal bundle of~\(Y\) in~\(X\) and let
  \(D_\nu\selfmap\colon \nu\to\nu\) be induced by the
  derivative of~\(\selfmap\) as above.  Let \(r_Y\colon
  \CONT_0(X)\to\CONT_0(Y)\) be the restriction map and let
  \(\Eul_Y\in\KK^G_0(\CONT_0(Y),\UNIT)\) be the equivariant
  Euler characteristic of~\(Y\).  Then
  \[
  \Lef(\selfmap) =
  r_Y \otimes_{\CONT_0(Y)} \sign(\ID_\nu- D_\nu\selfmap)
  \otimes_{\CONT_0(Y)} \Eul_Y.
  \]
  Furthermore, \(\Eul_Y\) is the equivariant
  \(\K\)\nb-homology class of the \emph{de Rham} operator
  on~\(Y\).
\end{theorem}

In brief, the Lefschetz invariant of \(G\) is the Euler characteristic 
of the fixed-point set, twisted by an appropriate equivariant line bundle 
depending on orientation data.

If \(\selfmap\) and~\(\ID_X\) are transverse then the 
fixed point subset~\(Y\) is discrete. A discrete set is a
manifold and its Euler characteristic -- a degenerate case 
of Proposition \ref{pro:euler} -- is represented by the 
Kasparov cycle in which the Hilbert space is 
 \(L^2(\Lambda^*_\C(T^*Y)) \defeq \ell^2(Y)\) equipped with the representation
\(\CONT_0(Y)\to\Comp(\ell^2Y)\) by pointwise multiplication
operators, and the zero operator. The group permutes the points of 
\(Y\) and so acts by unitaries on \(\ell^2(Y)\). 

The normal bundle~\(\nu\) to~\(Y\) in~\(X\) in this case is the restriction
of the vector bundle~\(\Tvert X\) to the subset~\(Y\).
  For \(p\in Y\), let~\(n_p\) be~\(+1\) if
\(\ID_{T_pX}-D_p\selfmap\) preserves orientation, and~\(-1\)
otherwise.  The graded equivariant line bundle
\(\sign(\ID_\nu-D_\nu\selfmap)\) in
Theorem~\ref{the:Lef_smooth} is determined by pairs
\((n_p,\chi_p)\) for \(p\in Y\), where~\(n_p\) is the parity of
the representation at~\(p\) and~\(\chi_p\) is a certain
real-valued character \(\chi_p\colon \Stab(p)\to\{-1,+1\}\)
that depends on \(\ID_{\Tvert_pX}-D_p\selfmap\) and the
representation of the stabiliser \(\Stab(p)\subseteq G\) on
\(\Tvert_pX\).  Equivariance
implies that~\(n_p\) is constant along \(G\)\nb-orbits,
whereas~\(\chi_p\) behaves like \(\chi_{g\cdot p} = \chi_p\circ
\operatorname{Ad}(g^{-1})\).  Let \(\ell^2_\chi(Gp)\) be the
representation of the cross-product \(G\ltimes\CONT_0\bigl(G/\Stab(p)\bigr)\)
obtained by inducing the representation~\(\chi_p\) from
\(\Stab(p)\), and let \(\CONT_0(X)\) act on \(\ell^2_\chi(Gp)\)
by restriction to \(G/\Stab(p)\).  This defines a
\(G\)\nb-equivariant \Star{}homomorphism
\[
\xi_{Gp,\chi}\colon \CONT_0(X)\to \Comp(\ell^2_\chi G).
\]
Theorem~\ref{the:Lef_smooth} asserts the following:

\begin{corollary}
\label{cor:maps_transverse}
  If the graph of~\(\varphi\) is transverse to the diagonal in
  \(X\times X\) then,
  \[
  \Lef(\selfmap) = \sum_{Gp\in G\backslash\Fix(\selfmap)}
  n_p [\xi_{Gp,\chi}]
  \]
  where \([\xi_{Gp,\chi}] \in \KK^G_0(\CONT_0(X),\UNIT)\) and
  the multiplicities~\(n_p\) are explained above.

  Furthermore, the character \(\chi\colon \Stab_G(p) \to
  \{-1,+1\}\) at a fixed point~\(p\) has the explicit formula
  \[
  \chi(g) = \sign \det \bigl(\ID - D_p\varphi\bigr)
  \cdot \sign \det \bigl( \ID - D_p\varphi_{\Fix(g)}\bigr).
  \]
\end{corollary}

If, in addition, \(G\) is trivial and~\(X\) is connected, then
\(\xi_{Gp,\chi}=\ev_p\) for all \(p\in Y\); moreover, all point
evaluations have the same \(\K\)\nb-homology class because
they are homotopic.  Hence we get the classical Lefschetz data
multiplied by the \(\K\)\nb-homology class of a point
\[
\Lef(\selfmap) =
\biggl(\sum_{p\in\Fix(\selfmap)}
\sign(\ID_{T_pX} - D\selfmap_p)\biggr) \cdot [\ev]
\]
as asserted above.  This sum is finite if~\(X\) is compact.

We include the following for the benefit of the reader; it 
enables her or she to verify our computations by direct inspection.

\begin{lemma}
 \label{lem:pointevals}
  Let \(H\subseteq G\) be compact and open, let \(p,q\in X^H\)
  belong to the same path component of the fixed point
  subspace~\(X^H\), and let \(\chi\in\Rep(H)\).  Then
  \[
  \bigl[\xi_{Gp,\ind_H^{\Stab(Gp)}(\chi)}\bigr] =
  \bigl[\xi_{Gq,\ind_H^{\Stab(Gq)}(\chi)}\bigr]
  \qquad \text{in \(\KK_0(\CONT_0(X),\UNIT)\).}
  \]
\end{lemma}

\begin{remark}
\label{rem:euler_is_points}
If the identity map \(\ID\colon X \to X\) can be equivariantly  
perturbed to be in general position in the sense explained above, 
then combining Proposition \ref{pro:euler} and Corollary 
\ref{cor:maps_transverse} proves that the class of the de Rham operator
in \(\KK^G(\CONT_0(X), \C)\) is a sum of 
point-evaluation classes. L\"uck and Rosenberg show in 
\cite{Lueck-Rosenberg:Lefschetz} that this can always be achieved when 
\(G\) is discrete and acts properly on \(X\).  
\end{remark}

The next two explicit examples involve isolated fixed points.

\begin{example}
  Let \(G\cong \Z\ltimes \Z/2\Z\) be the infinite dihedral
  group, identified with the group of affine transformations
  of~\(\R\) generated by \(u(x) = -x\) and \(w(x) = x+1\).
  Then~\(G\) has exactly two conjugacy classes of finite
  subgroups, each isomorphic to~\(\Z/2\).  Its action on~\(\R\)
  is proper, and the closed interval \([0,\nicefrac{1}{2}]\) is
  a fundamental domain.  There are two orbits of fixed point
  in~\(\R\), those of \(0\) and~\(\nicefrac{1}{2}\), and their
  stabilisers represent the two conjugacy classes of finite
  subgroups.

  Now we use some notation from
  Example~\ref{exa:orientation_representation}.  Each copy of
  \(\Z/2\) acting on the tangent space at the fixed point acts
  by multiplication by~\(-1\) on tangent vectors.  Therefore,
  the computations in
  Example~\ref{exa:orientation_representation} show that for
  any nonzero real number~\(A\), viewed as a linear
  transformation of the tangent space that commutes with
  \(\Z/2\), we have
  \[
  \sign (A) =
  \begin{cases}
    \R^\op_\chi&\text{if \(A <0\), and}\\
    \R_\tau&\text{if \(A>0\).}
  \end{cases}
  \]

  Let~\(\varphi\) be a small \(G\)\nb-equivariant perturbation of
  the identity map \(\R\to\R\) with the following properties.
  First, \(\varphi\) maps the interval \([0,\nicefrac{1}{2}]\) to
  itself.  Secondly, its fixed points in
  \([0,\nicefrac{1}{2}]\) are the end points \(0\),
  \(\nicefrac{1}{2}\), and \(\nicefrac{1}{4}\); thirdly, its
  derivative is bigger than~\(1\) at both endpoints and between
  \(0\) and~\(1\) at \(\nicefrac{1}{4}\).  Such a map
  clearly exists.  Furthermore, it is homotopic to the identity
  map, so that \(\Lef(\varphi)=\Eul_\R\).

  By construction, there are three fixed points modulo~\(G\),
  namely, the orbits of \(0\), \(\nicefrac{1}{4}\) and
  \(\nicefrac{1}{2}\).  The isotropy groups of the first and
  third orbit are non-conjugate subgroups isomorphic
  to~\(\Z/2\); from
  Example~\ref{exa:orientation_representation}, each of them
  contributes \(\R^\op_\chi\).  The point \(\nicefrac{1}{4} \)
  contributes the trivial character of the trivial subgroup.
  Hence
  \[
  \Lef (\varphi) = - [\xi_{\Z,\chi}] - [\xi_{\Z+\nicefrac{1}{2},\chi}]
  + [\xi_{\Z+\nicefrac{1}{4}}].
  \]

  On the other hand, suppose we change~\(\varphi\)
  to fix the same points but to have zero derivative at~\(0\)
  and~\(\nicefrac{1}{2}\) and large derivative at
  \(\nicefrac{1}{4}\).  This is obviously possible.  Then we
  get contributions of~\(\R_\tau\) at~\(0\) and
  \(\nicefrac{1}{2}\) and a contribution of
  \(-[\xi_{\nicefrac{1}{4}}]\) at \(\nicefrac{1}{4}\).  Hence
  \[
  \Lef (\varphi) = [\xi_{\Z, \tau}] +
  [\xi_{\Z+\nicefrac{1}{2},\tau}] - [\xi_{\Z+\nicefrac{1}{4}}].
  \]
  Combining both formulas yields the identity 
  \begin{equation}
    \label{mysteriousequality}
    [\xi_{\Z,\tau}] + [\xi_{\Z+\nicefrac{1}{2},\tau}]
    -[\xi_{\Z+\nicefrac{1}{4}}]
    = - [\xi_{\Z,\chi}] - [\xi_{\Z+\nicefrac{1}{2},\chi}]
    + [\xi_{\Z +\nicefrac{1}{4}}].
  \end{equation}
  By the way, the left-hand side is the description of
  \(\Eul_\R\) we get from the combinatorial dual with the
  obvious \(G\)\nb-invariant triangulation of~\(\R\) with
  vertex set \(\Z+\nicefrac{1}{2}\Z \subset \R\).

  Using Lemma \ref{lem:pointevals} one can check 
  \eqref{mysteriousequality} by direct computation. 
\end{example}

\section{Geometric KK-theory}
\label{sec:Geometric_KK}

In the previous chapter we explained our computation
 of the Lefschetz map for self-\emph{maps}
of a \(\Grd\)-space \(\Tot\) in several relatively simple situations. In 
these cases, \(\Grd\) was in each case a group (a groupoid with 
trivial base). Obviously, not all equivariant Kasparov 
self-morphisms \(\KK_*^\Grd\bigl( \CONT_0(\Tot), \CONT_0(\Tot)\bigr)\)
 are represented by maps.  We have organized this survey around the 
 problem of computing Lefschetz invariants of more general equivariant 
 Kasparov self-morphisms. This requires describing the morphisms 
 themselves in some geometric way. 
 The theory of \emph{correspondences} 
 of Baum, Connes and Skandalis (see \cite{Baum-Block:Bicycles} 
 \cite{Connes-Skandalis:Longitudinal}) would seem ideal for this 
 purpose. Since we are working in the equivariant setting, to use 
 it would necessitate checking that the pseudodifferential calculus which 
 plays such a prominent role in \cite{Connes-Skandalis:Longitudinal}
 works equivariantly with respect to a 
 an action of a group (or groupoid), as well as proving the 
 main functoriality result for \(\K\)-oriented maps.
  Although it seems plausible that such an extension could be carried out, 
  a major problem arises in connnection with composing 
 correspondences using transversality in the equivariant situation
 (we explain this below.) A trick of Baum and Block (see 
 \S\ref{subsec:top_category} and \cite{Baum-Block:Bicycles})
  is useful in this connection, but in order for it to work, some 
  hypotheses on vector bundles are necessary. Our approach is 
  to build in the vector bundle requirements into the definitions. 
  This is only reasonable if \(\Grd\) is proper; we now show how to 
  reduce to this case using the Baum-Connes conjecture. 
  
   \subsection{Using Baum-Connes to reduce to proper groupoids}
   \label{subsec:inflation_trick}
  Let \(G\) be a locally compact group (or groupoid). The \emph{classifying 
  space \(\EG{G}\) for proper actions of \(G\)} is the proper \(G\)-space with the 
  universal property that if \(X\) is any proper \(G\)-space, then there is a 
  \(G\)-equivariant classifying map \(\chi\colon X \to \EG{G}\) which is
   unique up to \(G\)-homotopy. If \(G\) is a proper groupoid to begin 
   with, then \(\EG{G}= \Base\) gives a simple model for 
   \(\EG{G}\), for it is proper as a \(G\)-space and has the required 
   universal property. In particular, if \(G\) is a compact group, 
   then \(\EG{G}\) is a point.
  
  For \(G\) and \(G\)-spaces \(X\) and \(Y\), the \emph{inflation map} 
  \eqref{eq:inflation_map}  
  \begin{multline}
  \label{mult:inflation_map}
   p_{\EG{G}}^*\colon \KK^G_*(\CONT_0(X), \CONT_0(Y)) \to 
  \RKK^G_*(\EG{G}; \CONT_0(X), \CONT_0(Y)) \\ \defeq 
  \KK^{G\ltimes \EG{G}}_*(\CONT_0(X\times \EG{G}), \CONT_0(Y\times \EG{G}))
  \end{multline}
  is an isomorphism as soon as the \(G\) action on \(X\) is topologically 
  amenable, and in particular as soon as it is proper. 
  An abstract (respectively Kasparov) dual for the \(G\)-space
   \(X\) pulls back to one for \(X\times \EG{G}\) as a 
   \(G\ltimes \EG{G}\)-space, and the diagram 
  \begin{equation}
  \label{eq:inflate_lef}
   \xymatrix{
 \KK^G_* ( \CONT_0(X), \CONT_0(X)) \ar[d]_{\cong}^{p_{\EG{G}}^*}
  \ar[r]^{\Lef} & \KK^G_*(\CONT_0(X), \C ) \ar[d]_{\cong}^{p_{\EG{G}}^*}\\
  \KK^{G \ltimes \EG{G}}_* \bigl( \CONT_0(X\times \EG{G}), 
  \CONT_0(X\times \EG{G}) )\ar[r]^{\Lef} & \KK^{G \ltimes \EG{G}}_*
  \bigl( \CONT_0(X\times \EG{G}), \CONT_0(\EG{G})\bigr)}
  \end{equation}
commutes. Hence the Lefschetz map for \(G\) acting on \(X\) is isomorphic to the 
Lefschetz map of \(G\ltimes \EG{G}\) acting on \(X\times \EG{G}\). 

This replaces the non-proper groupoid \(G\) by the proper groupoid 
\(\Grd\defeq G\times \EG{G}\) at no loss of information. 
  
In terms of this situation, our definitions are going to yield a 
a theory of \(G\ltimes \EG{G}\)-equivariant 
correspondences based on \(\K\)-oriented 
\(G\ltimes \EG{G}\)-equivariant vector bundles (equivalently, 
\(\Grd\)-equivariant vector bundles on \(X\times \EG{G}\))
and \(G\ltimes \EG{G}\)-equivariant 
open embeddings.  Such correspondences will yield analytic 
Kasparov morphisms since open 
embeddings do, while zero sections and projections 
of \(G\ltimes \EG{G}\)-equivariantly 
\(\K\)-oriented vector bundles yield analytic Kasparov morphisms 
for the Thom isomorphism for \(\K\)-oriented vector bundles over a space 
with an action of a proper groupoid, which is proved in \cite{LeGall:KK_groupoid}).
 
 In order to compose correspondences we need 
a sufficient supply of `trivial' vector bundles, but for this too, the fact 
that we have a proper groupoid makes a big difference. 
If \(\Grd\) is a groupoid acting 
on a space, then a \(\Grd\)-equivariant vector bundle (see 
\S \ref{subsec:equivariant_vector_bundles}) over that space is 
\emph{trivial} if it is pulled back from the unit space of \(\Grd\) using the anchor
map for the space. In the case 
we are discussing, where \(\Grd \defeq G\ltimes \EG{G}\),  
a \(G\ltimes \EG{G}\)-vector bundle on \(X\times \EG{G}\) is 
the same as a \(G\)-vector bundle on \(X\times \EG{G}\), and a 
trivial \(G\ltimes \EG{G}\)-vector bundle is a \(G\)-vector bundle 
on \(X\times \EG{G}\) which is pulled back from \(\EG{G}\) under 
the coordinate projection. 

In general, if a groupoid is not proper, 
it may have no equivariant vector bundles on its base, \emph{e.g.} if 
our groupoid is a group \(G\), then its base is a point, so that a trivial 
\(G\)-vector bundle over \(X\) is equivalent to a finite-dimensional representation 
of \(G\).

\subsection{Equivariant vector bundles}
\label{subsec:equivariant_vector_bundles}
 As per above, \(\Grd\) shall 
be a \emph{proper} groupoid until further notice. If \(\Grd\) is a 
group, this means that it must be compact. We will frequently 
consider this case in examples.

If \(X\) is a 
 \(\Grd\)-space, we remind the reader that 
 the anchor map for the action is denoted 
 \( \anchor_X\colon X \to \Base\). A 
\emph{\(\Grd\)-equivariant vector bundle on \(X\)} is a vector bundle on \(X\) which 
is also a \(\Grd\)-space such that elements of \(\Grd\) map fibres to fibres linearly.
There is an obvious notion of isomorphic \(\Grd\)-equivariant vector bundles. 

If \(\Grd\) is a compact group, then a \(\Grd\)-equivariant vector bundle on a 
point is a finite-dimensional linear representation of \(\Grd\). 

\begin{notation}
\label{notation}
If \(\VB\) is an equivariant vector bundle over \(X\) then we 
denote by \(\proj{\VB}\colon \VB \to X\) the vector bundle projection and 
\(\zers{\VB}\colon X \to \VB\) the zero section. We frequently denote by 
\(\total{\VB}\) the total space of \(\VB\), and denote by \(\vbK_\Grd (X)\) 
the Grothendieck group of the 
monoid of isomorphism classes of \(\Grd\)-equivariant vector bundles over 
\(X\). 
\end{notation}

Given that \(\Grd\) is assumed proper, 
equivariant vector bundles behave in some ways just like ordinary vector bundles. 
For example, if \(Y\subset X\) is a closed, \(\Grd\)-invariant subset of a 
\(\Grd\)-space \(X\), and if \(\VB\) is a \(\Grd\)-equivariant 
vector bundle on \(X\), then any equivariant section of \(\VB\) defined on \(Y\) 
can be extended to an equivariant section defined on
 an open \(\Grd\)-invariant neighbourhood of \(Y\). This involves an averaging 
 procedure (see \cite{Emerson-Meyer:Geometric_KK}).  As a consequence, 
 if \(f_i\colon X \to Y\), \(i=0,1\) are two \(\Grd\)-equivariantly homotopic maps, and if 
 \(\VB\) is a \(\Grd\)-equivariant vector bundle on \(Y\), then \(f_0^*(\VB)\) and 
 \(f_1^*(\VB)\) are \(\Grd\)-equivariantly isomorphic. 
 
 On the other hand, some new problems appear in connection with equivariant 
 vector bundles. We first formalize our notion of triviality and the corresponding 
 notion of \emph{subtriviality} in a definition. 
 
 \begin{definition}
 \label{def:trivial_subtrivial}
 Let \(X\) be a \(\Grd\)-space. 
 A \(\Grd\)-vector bundle over \(X\) is \emph{trivial} if it is pulled back from \(\Base\)
 under the anchor map \(\anchor_X\colon X \to \Base\). We denote the pull-back of a 
 \(\Grd\)-vector bundle \(E\) over \(\Base\) to the \(\Grd\)-space \(X\) by \(E^X\).  A
 \(\Grd\)-equivariant vector bundle is \emph{subtrivial} if it is a direct summand of a 
 trivial \(\Grd\)-vector bundle. 
 
 \end{definition}
  
  \begin{example}
  \label{ex:trivial_bundles}
 If \(\Grd\) is a compact group, then a trivial \(\Grd\)-equivariant 
 vector bundle over \(X\) has the form \(X\times \R^n\) where \(\R^n\) carries 
 a linear representation of \(\Grd\), and where \(\Grd\) acts on \(X\times \R^n\)
 diagonally.  
 
 \end{example}

We will make repeated use of the following basic fact about representations of 
 compact groups (see \cite{Mostow:Equivariant_embeddings}). 
  
 \begin{lemma}
 \label{lem:basic_fact}
 Let \(\Grd\) be a compact group and \(\Grd' \subset\Grd\) be a subgroup. Then 
 any finite-dimensional representation of \(\Grd'\) is contained in 
 the restriction to \(\Grd'\) of a finite-dimensional representation of \(\Grd\). 
 \end{lemma}

Even if \(\Grd\) is a compact group, 
due to the notion of `trivial' vector bundle we are using, not every \(\Grd\)-vector bundle 
is locally trivial in the category of \(\Grd\)-vector bundles, \emph{i.e.} locally 
isomorphic 
to a trivial \(\Grd\)-vector bundle.  
But it is not hard to prove the following and it is a good exercise for 
understanding equivariant vector bundles. We leave the proof to 
the reader, but see \S \ref{subsec:intro_to_embedding_theorems} for 
some important ingredients in the argument. 
 
\begin{lemma}
\label{lem:locally_subtrivial}
Every \(\Grd\)-equivariant vector bundle on \(X\) is locally subtrivial in the sense that 
for every \(x \in X\) there is a \(\Grd\)-equivariant vector bundle on \(\Base\), a 
\(\Grd\)-equivariant neighbourhood \(U\) of \(x\), and an embedding
\(\varphi\colon \VB_{|_{U}}\to E^U\). 
\end{lemma}

Improving this local result 
to a global one is not possible, however, without an appropriate 
compactness assumption.

\begin{example}
\label{ex:not_subtrivial}
Let \(X\defeq \Z\) with the trivial action of the compact group 
\(\Grd\defeq \mathbb{T}\). Then the \(1\)-dimensional complex vector 
bundle \(\Z\times \C\) with the action of \(z\in \Grd \) in the fibre over 
\(n\) by the character \(z\mapsto z^n\) is not subtrivial, since it contains infinitely 
many distinct irreducible representations of \(\Grd\). 

\end{example}

\begin{definition}
\label{def:enough_vbs}
Let \(X\) be a \(\Grd\)-space. 
\begin{itemize}
 \item The space \(X\) \emph{has enough \(\Grd\)-equivariant vector bundles} if whenever
 we are given \(x\in X\) and a finite-dimensional representation of the compact isotropy group 
 \(\Grd^x_x\) of \(x\), there is a \(\Grd\)-equivariant vector bundle over \(X\) whose restriction 
 to \(x\) contains the given representation of \(\Grd^x_x\).  
 \item The space \(X\) has a \emph{full} vector bundle if there is a vector bundle 
 \(\VB\) over \(X\) such that any irreducible representation of \(\Grd^x_x\) is contained 
 in the representation of \(\Grd^x_x\) on \(\VB_x\) (and we say that such a vector bundle 
 \(\VB\) is full.)
 \end{itemize}
\end{definition}

It is the content of Lemma \ref{lem:basic_fact} that a compact group acting on a 
point has enough vector bundles.  It does not have a full vector bundle unless it is finite, 
because a compact group with a finite dual has finite-dimensional \(L^2(G)\) by the 
Peter-Weyl theorem and must then be finite.

It is easy to check that if 
 \(f\colon X \to Y\) is a \(\Grd\)-equivariant map 
then \(X\) has enough equivariant vector bundles if \(Y\) does, and 
\(f^*(\VB)\) is a full vector bundle on \(X\) if \(\VB\) is a full vector bundle on 
\(Y\). Both of these assertions use the basic fact Lemma \ref{lem:basic_fact}.

\begin{example}
\label{ex:no_full}
The \(\Grd\)-space described in Example \ref{ex:not_subtrivial} does not have a full 
vector bundle, although it obviously has enough 
vector bundles.
\end{example}

The following example is more subtle. It is 
due to Julianne Sauer (see \cite{Sauer:K-theory}).

\begin{example}
\label{ex:not_enough}
Let \(X = \R\) and \(K\) be the compact group \(K \defeq 
\prod_{n \in \Z} \Z/2\) acting \emph{trivially} on \(\R\). Since \(X\) is 
\(K\)-equivariantly contractible, and by the 
homotopy-invariance of equivariant 
vector bundles (mentioned above at the beginning of
 \S \ref{subsec:equivariant_vector_bundles}), 
any \(K\)-equivariant vector bundle \(\VB\) on \(X\) is trivial, and hence 
all the
representations of \(K\) on the fibres \(\VB_x\) 
are equivalent.  

Now let \(\sigma \colon K \to K\) be the shift automorphism and consider
the group 
\(G \defeq K \ltimes_\sigma \Z\); it acts on \(X\) by letting \(\sigma (t) \defeq t+1\). 
This is a proper action.
Of course as a groupoid \(G\) is not proper, but we repair this below. 

We claim that the only trivial \(G\)-vector bundles on \(X\) yield trivial \(K\)-representations
in their fibres. This will show that \(X\) does not have enough 
\(G\)-vector bundles. 

The proof is as follows: any \(G\)-vector bundle \(\VB\)
on \(X\) must be trivial as a \(K\)-vector bundle, as above. 
On the 
other hand, the covariance rule for the semi-direct product implies that 
representations of \(K\) on \(\VB_x\) and \(\VB_{x+1}\) are mapped to each other 
(up to equivalence) by the action of \(\sigma\), and therefore 
\(\hat{\sigma}\colon \Rep (K) \to \Rep (K)\) fixes the point \([\VB_x]\).  
But \(\Rep (K)\) is the direct \emph{sum} of \(\widehat{\Z/2}\)'s and 
\(\widehat{\sigma}\) acts as the shift. The only fixed point then is 
the zero sequence. This corresponds to a trivial representation. 
This means that at every point \(x\in X\) the representation of 
\(K\) we get on \(\VB_x\) is trivial.

To repair the non-properness of \(G\), replace it by 
\(\Grd\defeq G\ltimes \EG{G}\) and replace \(X\) by \(X\times \EG{G}\) as
explained at the beginning of this section, then we get an example of a 
proper groupoid \(\Grd\) and a \(\Grd\)-space \(X\times \EG{G}\) such which 
does not have enough equivariant vector bundles. This is because 
\(X\times \EG{G}\) and \(X\) are \(G\)-equivariantly homotopy-equivalent 
anyway, because the action of \(G\) on \(X\) is proper. Hence these 
spaces have canonically isomorphic monoids of isomorphism classes of 
equivariant vector bundles. 

A Morita-equivalent approach is via a mapping cyclinder construction 
and produces a compact groupoid acting on a compact space without 
enough vector bundles. 
Take \([0,1]\times K\) modulo the relation \((1,k) \sim (0,\sigma (k))\). 
This results in a bundle of compact groups over the circle which can 
shown to be locally compact groupoid with Haar system.

 Let 
\(\Grd\) be this groupoid: it is proper. Its base \(\Base\) is the circle.
 By a holonomy argument similar to the one just given, any 
\(\Grd\)-equivariant vector bundle over \(\Base\) must restrict in each 
fibre to a trivial representation of \(K\). Thus, there are not enough 
\(\Grd\)-equivariant vector bundles on \(\Base\). 

\end{example}

 \begin{example}
\label{ex:discrete_enough}
If \(G\) is a discrete group with a \(G\)-compact model for 
\(\EG{G}\), then L\"uck and Oliver have shown in 
\cite{Lueck-Oliver:Completion} that 
 there is a full \(\Grd\)-equivariant vector bundle on 
\(\Base\), where
\(\Grd\defeq G \ltimes \EG{G}\), (so that \(\Base = \EG{G}\).)
\end{example}

\subsection{The topological index of Atiyah-Singer}
\label{subsec:intro_to_embedding_theorems}
We now indicate why the condition of having enough vector bundles, or having 
a full vector bundle, is important for describing analytic equivariant \(\KK\)-groups 
topologically. 

We start with the problem, famously treated by 
Atiyah and Singer, of describing the equivariant (analytic) 
index of a \(\Grd\)-equivariant 
elliptic operator in topological terms, where \(\Grd\) is a compact group, keeping in 
mind that an equivariant elliptic operator is an important example of a 
cycle for equivariant 
\(\KK\)-theory. 
 
Let \(X\) be a smooth manifold with a smooth action of the compact group 
\(\Grd\). The symbol 
of an equivariant elliptic operator on \(X\) is an equivariant \(\K\)-theory class 
for \(\Tvert X\). The idea of Atiyah and Singer 
for defining the topological index of the operator is to smoothly embed \(X\) 
in a finite-dimensional (linear) representation of \(\Grd\) on \(\R^n\). The
 derivative of this embedding gives a smooth embedding of \(\tang\) in  \(\R^{2n}\), 
 where \(\tang\) has the induced action of \(\Grd\).  
 Since \(\tang\) is an equivariantly \(\K\)-oriented manifold, the normal bundle \(\nu\)
 to the embedding is a \(\Grd\)-equivariantly \(\K\)-oriented vector bundle on 
 \(\tang\). The tubular neighbourhood embedding identities it with an open, 
 \(\Grd\)-equivariant neighbourhood of the image of \(\tang\) in \(\R^{2n}\). We 
 now obtain a composition 
 \[ \K^*_\Grd(\tang)\rightarrow \K^{*}_\Grd(N) \rightarrow 
 \K^*_\Grd(\R^{2n}) \to \K^*_\Grd(\star)\cong \Rep (\Grd),\]
 where the first map is the Thom isomorphism for the equivariantly \(\K\)-oriented 
 \(\Grd\)-vector bundle \(N\), the second is the map on equivariant \(\K\)-theory 
 induced by the open inclusion \(N \opem \R^{2n}\) and the third is 
 equivariant Bott Periodicity (\(\R^{2n} \cong \C^n\) 
 with the given action of 
 \(\Grd\) has an equivariant complex structure, so a 
 \(\Grd\)-equivariant \(\textup{spin}^c\)-structure. The spinor bundle 
 is the trivial \(\Grd\)-vector bundle \(\Lambda^*_\C (\C^n)\) over 
 \(\C^n\).) 
 
 The content of the index theorem is that this composition 
 agrees with the map \(\K^*_\Grd(\tang) \to \Rep (\Grd)
 \) obtained by first interpreting cycles for \(\K^*_\Grd(\tang) \) as symbols of equivariant 
 elliptic operators on \(X\), making these elliptic operators into Fredholm 
 operators, and taking their equivariant indices. 
 
But how do we get a smooth, equivariant
 embedding of \(X\) in a finite-dimensional linear representation of \(\Grd\) 
 in the first place? Since it involves important ideas for us, we will sketch the 
 proof. The result seems due to Mostow 
 (see \cite{Mostow:Equivariant_embeddings}). Very similar arguments also prove
  Lemma \ref{lem:locally_subtrivial}.

First of all, we may assume (or average using the Haar system on \(\Grd\))
 that \(X\) has an invariant 
Riemannian metric. Now the orbit \(\Grd x\) is a smooth embedded submanifold
of \(X\) isomorphic to \(\Grd/\Grd^x_x\). The tangent space of \(X\) at \(x\) splits into the 
orthogonal sum of the tangent space to the orbit and its orthogonal complements
\(N_x\defeq T_x(\Grd x)^\perp\). The latter is a finite-dimensional representation 
of \(\Grd^x_x\), and inducing it results in a \(\Grd\)-equivariant vector bundle
 \[N \defeq \Grd \times_{\Grd^x_x} N_x \defeq \Grd\times N_x \,/ \,(g,n) \sim 
 (gh, h^{-1}n)\; \text{for}\; h \in \Grd^x_x \]
on the orbit which is precisely the normal bundle to the embedded submanifold
\(\Grd x\).  By exponentiating we obtain an equivariant diffeomorphism between 
the total space of \(N\) and an invariant open neighbourhood \(U\) 
of the orbit. We embed this neighbourhood as follows.

By Lemma \ref{lem:basic_fact}, the representation 
of \(\Grd^x_x\) on \(N_x\) is contained in the restriction 
of some representation of \(\Grd\) on some finite-dimensional vector space \(\tilde{N}_x\).
The naturality of induction implies that we have an inclusion of vector bundles 
\(N \subset \tilde{N} \defeq \Grd\ltimes_{\Grd^x_x} \tilde{N}_x\). But 
since \(\tilde{N}_x\)  is the restriction of a \(\Grd\)-representation, 
\(\tilde{N}\) is a product bundle, \emph{i.e.} a trivial \(\Grd\) 
vector bundle on the orbit. This provides a \(\Grd\)-equivariant map  
\(U \to \tilde{N}_x\), explicitly, by mapping \([(g, n)]\in U \cong N \subset \tilde{N}\) to 
the point \(gn\in \tilde{N}_x\). 
It is of course not necessarily an embedding; to improve it to an 
embedding, fix a vector 
\(v\in \tilde{N}_x\) whose isotropy in \(\Grd\) is \emph{exactly} \(\Grd^x_x\) 
(for this see also \cite{Mostow:Equivariant_embeddings}) and set 
\[\varphi \colon U \cong \Grd\times_{\Grd^x_x}N_x\to \tilde{N}_x\oplus \tilde{N}_x, 
 \hspace{1cm} \varphi ([g,n]) \defeq (gn, gv).\]
The map \(\varphi\) is an equivariant embedding as required.

As mentioned, if \(X\) is compact, then we can then we can (carefully) 
paste together the local embeddings to get an embedding of \(X\); see 
\cite{Emerson-Meyer:Geometric_KK}, or the source
\cite{Mostow:Equivariant_embeddings}. 
 The reader should notice that 
the assumption that \(X\) has enough vector bundles 
is used implicitly to show that the representation of 
\(\Grd^x_x\) on \(N_x\) can be extended to a \(\Grd\)-equivariant vector bundle 
on the orbit of \(x\) (the vector bundle \(\tilde{N}\) induced from \(\tilde{N}_x\)). 
This was the statement of  
Lemma \ref{lem:basic_fact}, and is just an explicit way of saying that 
\(\Grd\) has enough vector bundles on its one-point base space.

 \subsection{Embedding theorems from \cite{Emerson-Meyer:Geometric_KK}}
 More generally, in \cite{Emerson-Meyer:Geometric_KK} the following is proved. 
 Let \(X\) be a \(\Grd\)-space, where \(\Grd\) is a proper groupoid.

 We say that \(X\) is a \emph{smooth \(\Grd\)-manifold} if 
 we can cover \(X\) by charts of the form 
 \(U\times \R^n\) where \(U\subset \Base\) is open, 
 so that with respect to this product structure the anchor map \(\anchor_X\colon X \to \Base\) 
 identifies with the first coordinate 
 projection, and such that groupoid elements and change of coordinates 
 are smooth in the vertical direction. 
 
 An \emph{smooth open embedding} between \(\Grd\)-manifolds is 
 a smooth equivariant map which is a diffeomorphism onto an open 
 subset of its codomain. 
  
\begin{theorem}
\label{thm:normal_factorisation}
Let \(\Grd\) be a (proper) groupoid and \(X\) and \(Y\) be 
smooth \(\Grd\)-manifolds. 

Suppose that either 
\begin{itemize}
\item[{\bf A}.] The \(\Grd\)-space
 \(\Base\) has enough vector bundles and \(\Grd\backslash X\) is compact, 
 \item[{\bf B}.] \(\Base\) has a full vector bundle and \(\Grd\backslash X\) has finite
covering dimension.  
\end{itemize}

Then, given a smooth, \(\Grd\)-equivariant map \(f\colon X \to Y\), there exists 
\begin{itemize}
\item A smooth \(\Grd\)-equivariant vector bundle \(\VB\) over \(X\),
\item A smooth \(\Grd\)-equivariant vector bundle \(\Triv\) over \(\Base\),  
\item An smooth, equivariant open embedding \(\varphi\colon \VB\to \Triv^Y\), 
\end{itemize}
 such that 
 \begin{equation}
 \label{eq:normal_factorisation}
  f = \proj{\Triv^Y}\circ \varphi \circ \zers{\VB}.
 \end{equation}
Furthermore, under the union of hypotheses \(\bf{A},\bf{B}\), any \(\Grd\)-equivariant 
vector bundle over \(X\) (or \(Y\)) is subtrivial. 
\end{theorem}

Recall that the notation \(\Triv^Y\) means the pullback of \(\Triv\) to 
\(Y\) using the anchor map \(\anchor_Y\colon Y \to \Base\).

We call a factorisation of a map \(f\colon X\to Y\) of the form 
\eqref{eq:normal_factorisation} a \emph{normal factorisation}.

\subsection{Normally non-singular maps}
\label{subsec:normal_maps}
As in the previous section, \(\Grd\) is a proper groupoid. The 
constructions of the previous section motivate the following definition. 

\begin{definition}
\label{def:normal_map}
A \emph{\(\Grd\)-equivariant normally non-singular map  \(\Phi\)} from \(X\) to \(Y\) is a 
triple \((\VB,  \Triv, \hat{f})\) where \(\VB\) is a \(\Grd\)-equivariant 
subtrivial vector bundle over \(X\), \(\Triv\) is a \(\Grd\)-equivariant 
vector bundle over \(\Base\) and \(\hat{f} \colon \VB\to \Triv^Y\) is 
a \(\Grd\)-equivariant open embedding. 
\begin{itemize}
\item The \emph{trace} of \(\Phi\) is the composition 
\(\proj{\Triv^Y}\circ \hat{f} \circ \zers{\VB}\).
\item The \emph{stable 
normal bundle} of \(\Phi\) is the class \([\VB] - [\Triv^X]\in \vbK_\Grd (X)\).
\item The \emph{degree} of \(\Phi\) is
\(\dim (\VB) - \dim (\Triv)\). 
\item The normally non-singular map  \(\Phi\) is \emph{\(\K\)-oriented}
 if the \(\VB\) and \(E\) are equivariantly \(\K\)-oriented. 
 \item The normally non-singular map  \(\Phi\) is \emph{smooth}
if \(X\) and \(Y\) are smooth \(\Grd\)-manifolds  
and \(\VB\) and \(\Triv\) are smooth equivariant vector bundles on which 
\(\Grd\) acts smoothly, and if \(\hat{f}\) is a smooth embedding. 
\end{itemize}
\end{definition}

Of course the trace of a \(\Grd\)-equivariant smooth normally non-singular map  is itself a 
smooth equivariant map, and the content of Theorem \ref{thm:normal_factorisation}
is that, conversely, any smooth equivariant map between smooth \(\Grd\)-manifolds 
is the trace of some smooth normally non-singular map , under some hypotheses about the 
availability of equivariant vector bundles. (This statement is improved in Theorem 
\ref{thm:normal_map_unique}.)

\begin{example}
\label{ex:simplest_normal}
The simplest example of a normally non-singular map  is the zero section and bundle projection 
\[ \zers{\VB}\colon \Tot \to \VB, \hspace{1cm} \proj{\VB}\colon \total{\VB} \to \Tot,\]
of a \(\Grd\)-equivariant vector bundle \(\VB\) over a \emph{compact} \(\Grd\)-space \(\Tot\), 
where \(\Grd\) is a compact group. The zero section is the trace of the normally non-singular map  
\( (\VB, 0_{\total{\VB}}, \ID)\). The stable normal bundle is the class 
\([\VB] \in \vbK_\Grd(X)\in \) of the vector bundle itself. Since \(X\) is compact, \(\VB\) is subtrivial. 
If \(X\) is not compact, this can fail, \emph{c.f.} Example \ref{ex:not_subtrivial}.

With the same (compact) \(X, \VB\) \emph{etc}, 
the bundle projection \(\proj{\VB}\colon \total{\VB} \to 
\Tot\) is the trace of a normally non-singular map   
\( (\pi_\VB^*(\VB'), \Triv, \varphi)\), where \(\VB'\in \Vect_\Grd(\Tot)\) is a choice of 
\(\Grd\)-vector bundle on \(\Tot\) such that \(\VB\oplus \VB' \) is a trivial bundle 
\(\Triv^\Tot\), and \(\varphi \colon \total{\VB\oplus \VB'}\cong \total{\proj{\VB}^*(\VB')}
\xrightarrow{\cong} \Triv^X\) is a trivialisation. 
The stable normal bundle is  
\( \proj{\VB}^*([\VB']) - [\Triv^\Tot] \in \vbK_\Grd(\total{\VB})\) respectively.

The normally non-singular map  just described seems to depend on the choice of trivialisation of 
\(\VB\), but it can be checked that any two 
choices yield equivalent normally non-singular maps 
in the sense explained below. 
\end{example}

\begin{example}
 \label{ex:normal_classical}
 Let \(\Grd\) be a compact group acting smoothly on 
 manifolds \(X,Y\) with \(X\) compact.
  By the discussion in \S \ref{subsec:intro_to_embedding_theorems} we can fix 
  a smooth, equivariant embedding \(i\colon X \to \Triv\) in a linear representation 
  of \(\Grd\). Define \(\VB\) to be the normal bundle to the embedding 
  \(x\mapsto \bigl( f(x), i(x)\bigr)\) of \(X\) in \(E^Y \defeq Y\times E\). 
  Let \(\varphi\colon \VB\to E^Y\) be the corresponding  
  tubular neighbourhood embedding. Then the trace of the composition 
  \(\proj{E^Y}\circ \varphi \circ \zers{\VB}\) of the normally non-singular map  
  \((\VB, \Triv, \hat{f})\) is \(f\). Since \(\tang \oplus \VB
  \cong f^*(\Tvert Y) \oplus E^X\), the stable normal bundle is 
  \(f^*([\Tvert Y])-[\tang]\in \vbK_\Grd(X)\). 
 
 \end{example}

\begin{example}
If \(G\) is a discrete group with a \(G\)-compact model for 
\(\EG{G}\), then L\"uck and Oliver have shown in 
\cite{Lueck-Oliver:Completion} that 
 there is a full \(\Grd\)-equivariant vector bundle on 
\(\EG{G}\), where
\(\Grd\defeq G \ltimes \EG{G}\), (so that the base of \(\Grd\) is 
\(\Base \defeq \EG{G}\).)
Let \(X\) and \(Y\) be smooth manifolds equipped with smooth actions of \(G\)
and \(f\colon X \to Y\) be a smooth, \(G\)-equivariant map. 
As above let \(\Grd\) be the proper groupoid \(G\times \EG{G}\). 
Applying the Baum-Connes procedure of \S \ref{subsec:inflation_trick} 
to this situation we get smooth 
 \(\Grd\)-manifolds \(X\times \EG{G}\) and 
\(Y\times \EG{G}\) and a smooth \(\Grd\)-map 
\(f\times \ID_{\EG{G}}\colon 
X\times \EG{G} \to Y\times \EG{G}\).
It is the trace of a normally non-singular map  because of 
Theorem \ref{thm:normal_factorisation} and the result of L\"uck and Oliver. 

 As we will see in the next section, if \(f\) is also \(\K\)-orientable in an appropriate 
 sense, then it will give rise to a morphism in \(\KK^\Grd(\CONT_0(\EG{G}\times X), 
 \CONT_0(\EG{G}\times Y).\) If \(X\) is a topologically amenable 
 \(G\)-space, this gives an element of 
 \(\KK^G(\CONT_0(X), \CONT_0(Y))\). 
\end{example}

Two normally non-singular map s are \emph{isomorphic} if there are vector bundle isomorphisms 
\(\VB_0\cong \VB_1\) and \(\Triv_0\cong \Triv_1\) that intertwine the open embeddings
\(f_0\) and \(f_1\). The \emph{lifting} of a normally non-singular map  
\(\Phi = \Psi = (\VB, \varphi, \Triv)\) along an equivariant vector bundle 
\(\Triv^+\) over \(\Base\) is the normally non-singular map  
\(\Phi \oplus \Triv^+ \defeq 
(\VB\oplus (\Triv^+)^X, \Triv\oplus \Triv^+, \hat{f}\times_\Base \ID_{\Triv^+})\). 
Two normally non-singular map s are \emph{stably isomorphic} if there are \(\Grd\)-equivariant 
vector bundles \(\Triv_0^+\) and \(\Triv_1^+\) such that \(\Phi_0\oplus \Triv^+_0\) is 
isomorphic to \(\Phi_1\oplus \Triv^+_1\). Finally, two 
normally non-singular map s \(\Phi_0\) and \(\Phi_1\) 
are \emph{isotopic} if there is a continuous 
\(1\)-parameter family of normally non-singular map s whose values at the endpoints are 
stably isomorphic to \(\Phi_0\) and \(\Phi_1\) respectively 
(see \cite{Emerson-Meyer:Geometric_KK} for the exact definition), and are 
\emph{equivalent} if they have isotopic liftings. There is an obvious notion 
of \emph{smooth equivalence} of smooth normally non-singular map s. 

There are obvious \(\K\)-oriented analogues of the above relations. For 
example, lifting must only use \(\K\)-oriented trivial bundles, and isomorphism 
must preserve the given \(\K\)-orientations. Referring to this 
kind of equivalence we will speak of \emph{\(\K\)-oriented equivalence} of 
\(\K\)-oriented normally non-singular maps.

\subsection{Manifolds with smooth normally non-singular maps to \(\Base\)}
\label{subsec:normal_manifolds}
A useful hypothesis covering a number of geometric situations is that 
a given smooth \(\Grd\)-manifold \(X\) admits a smooth normally non-singular map  to 
the object space \(\Base\) of \(\Grd\). By the theorem above this
 is the case if \({\bf A}\) or \(\bf{B}\) hold. It means explicitly that 
we have a triple \( (\Normal_X, \hat{g}, \Triv)\) where \(\Normal_X\) is a 
smooth subtrivial vector bundle over \(X\), \(\Triv\)
 is an equivariant vector bundle 
over \(\Base\) and \(\hat{g}\) is a smooth open equivariant 
embedding \(\Normal_X \to \Triv\). 
Note that \(\Normal_X\oplus \tang
\cong\Triv^X\). Such a normally non-singular map  is (smoothly) stably isomorphic 
to a \(\K\)-oriented normally non-singular map  because we can replace 
if needed \(\Triv\) by \(\Triv\oplus \Triv\), which is canonically 
equivariantly \(\K\)-oriented using the \(\Grd\)-equivariant complex structure, 
and replacing \(\Normal_X\) by \(\Normal_X\oplus \Triv^X\). 

If \( ( \Normal_X, \hat{g}, \Triv) \) is a smooth normal
 map to \(\Base\) such that \(\Triv\) is equivariantly \(\K\)-oriented, then 
 \(\K\)-orientations on \(\Normal_X\) are in 1-1 correspondence 
 with \(\K\)-orientations on \(\tang\) because of the 2-out-of-3 property. 
 One can prove the following.

 \begin{theorem}
  \label{thm:normal_map_unique}
  Let \(\Source\) and~\(\Target\) be smooth
  \(\Grd\)\nb-manifolds, and assume that \(X\) admits a smooth, normal 
  \(\Grd\)-map to \(\Base\) and that \(f^*(\Tvert Y)\) is subtrivial. 
  
  Then any smooth \(\Grd\)\nb-map
  from~\(\Source\) to~\(\Target\) is the trace of a smooth
  normal \(\Grd\)\nb-map, and two smooth normally non-singular maps
  from~\(\Source\) to~\(\Target\) are smoothly equivalent if
  and only if their traces are smoothly homotopic.
  
  Furthermore, smooth equivalence classes of 
  smooth \(\K\)-oriented normally non-singular maps 
  from \(X\) to \(Y\) are in 1-1 correspondence with pairs \( (f, \tau)\) 
  where \(f\) is a smooth homotopy class of equivariant smooth map 
  \( X\to Y\) and \(\tau \) is an equivariant \(\K\)-orientation on 
  \(\Normal_X \oplus f^*(\tang)\). 
  
\end{theorem}
 
 We sketch the existence part of this proof. Fix
 a smooth equivariant normally non-singular map 
  \( (\Normal_X, \hat{g}, \Triv)\) from 
 \(X\) to \(\Base\). We can assume by replacing \(\Triv\) by 
 \(\Triv \oplus \Triv\) and \(\Normal_X\) by \(\Normal_X\oplus \Triv^X\) 
 if needed that \(\Triv\) is equivariantly \(\K\)-oriented. Let
   \(Y\) be another smooth \(\Grd\)-manifold and \(f\colon X \to Y\) 
 be a smooth map. 
  Let \(g\colon X \xrightarrow{\zers{\Normal_X}} \Normal_X\xrightarrow{\hat{g}} E\) the 
 composite smooth embedding. 
 
 One obtains a 
 a smooth embedding
  \(X \to  Y\times_\Base \Triv =  \Triv^Y\), \(x\mapsto \bigl( f(x), g(x)\bigr)\). 
  It has a (smooth) normal bundle \(\VB\) with a smooth
   open embedding in \(\Triv^Y\). 
  Since \(V\cong \Normal_X\oplus f^*(\Tvert Y)\), \(\VB\) is subtrivial and 
  \((\VB, \Triv, \hat{f})\) is a smooth normally non-singular map  with trace \(f\) and 
  stable normal bundle 
  \( [\VB]-[\Triv^X]\in \vbK_\Grd (X)\). Note that \(\VB\oplus \tang 
  \cong f^*(\Tvert Y)\oplus \Triv^X\). 
  
  The stable normal bundle is \( [\VB] -[\Triv^X] = f^*([\Tvert Y]) - [\tang]\in 
  \vbK_\Grd(X)\). Equivariant 
  \(\K\)-orientations on \(\VB\) are in 1-1 correspondence with equivariant 
  \(\K\)-orientations on \(\Normal_X\oplus f^*(\Tvert Y)\).

\subsection{Correspondences}
We are now in a position to define what correspondences are. 
Let \(\Grd\) continue to denote a proper groupoid. 

\begin{definition}
\label{def:correspondence}
Let \(X\) and \(Y\) be \(\Grd\)-spaces. 
A \emph{\(\Grd\)-equivariant correspondence from \(X\) to \(Y\)}
 is a quadruple 
\(( M, \mapl, \mapr, \xi)\) where \(M\) is a \(\Grd\)-space, 
\(\mapr\colon M \to Y\) is a \(\Grd\)-equivariantly \(\K\)-oriented 
normally non-singular map , \(\mapl\colon M \to X\) is an equivariant map, and 
\(\xi \in \RK^*_{\Grd, X}(M)\) is a \(\Grd\)-equivariant \(\K\)-theory 
class with \(X\)-compact support (see \S \ref{subsec:equivariant_K-theory})
 where the \(\Grd\ltimes X\)-structure on \(M\)
is that determined by the \(\Grd\)-equivariant map \(\mapl\colon M \to X\). 

The \emph{degree } of the correspondence \(( M, \mapl, \mapr, \xi)\)
is the sum of the degrees of \(\mapr\) and \(\xi\). 

\end{definition}

\begin{remark}
Thus a significant difference from the set-up of Connes and Skandalis 
in \cite{Connes-Skandalis:Longitudinal} is that the map 
\(\mapl\colon M \to X\) is not required to be proper; we have replaced 
this by a support condition on \(\xi\).
\end{remark}

Several equivalence relations on correspondences are imposed. 
The first is 
to consider two correspondences
\((M, \mapl_0, \mapr_0, \xi)\) and \((M, \mapl_1, \mapr_1, \xi)\)  
to be equivalent if their normally non-singular maps are equivalent. 
The second is to consider \emph{bordant} correspondences equivalent (we 
will not discuss this at all in this survey.) 
The third is most interesting, and is called 
\emph{Thom modification}. The Thom modification of a 
correspondence \((M, \mapl, \mapr, \xi)\) using a subtrivial \(\K\)-oriented 
vector bundle \(\VB\) over \(M\) is the correspondence 
\[ \bigl( \VB, \mapl\circ \proj{\VB}, \mapr\circ \proj{\VB}, \tau_\VB (\xi)\bigr),\]
where 
\(\tau_\VB\colon \RK^*_{\Grd, X}(M) \xrightarrow{\cong}
 \RK^{* +\dim(\VB)}_{\Grd, X}(\total{\VB})\)
is the Thom isomorphism, \emph{i.e.} \(\tau_\VB (\xi) \defeq \proj{\VB}^*(\xi)\cdot \xi_\VB\), 
where \(\xi_\VB\in \RK^{\dim (\VB)}_{\Grd, M}(\total{\VB})\) is the Thom class. 
We declare a correspondence and its Thom modification to be Thom 
equivalent. Note that applying Thom modification to a correspondence 
does not change its degree. 

The equivalence relation on correspondences is that generated by equivalence 
of normally non-singular maps, bordism and Thom equivalence. 

\begin{definition}
\label{def:topological_KK}
Let \(\Grd\) be a proper groupoid and \(X\) and \(Y\) be \(\Grd\)-spaces. 
We let \(\GKK_\Grd^*(X,Y)\) denote the \(\Z/2\)-graded set of equivalence classes of 
\(\Grd\)-equivariant correspondences from \(X\) to \(Y\), graded by degree. 
\end{definition}

A correspondence \( (M, \mapl, \mapr, \xi)\) from \(X\) to \(Y\) is \emph{smooth} if 
\(X, Y\) and \(M\) are smooth manifolds, \(\mapr\) is a smooth normally non-singular map 
(see \S \ref{subsec:intro_to_embedding_theorems}), and 
\(\mapl\) is a smooth map. There is a rather obvious notion of \emph{smooth equivalence}
of smooth correspondences. This gives rise to a parallel theory using only 
smooth equivalence classes of smooth correspondences; we do not use notation for 
this.

\subsection{\(\GKK^\Grd\) as a category}
\label{subsec:top_category}
Classes of correspondences form a category with analogous properties to 
Kasparov's equivariant \(\KK\) (that is, to \emph{analytic} Kasparov theory). 
The composition 
of correspondences is called the \emph{intersection product}. For composition we
use similar notation to Kasparov's: if 
\(\Psi \in \GKK^\Grd_i(X, Y)\) is a topological morphism, \emph{i.e.} an 
equivalence class of equivariant correspondence from \(X\) to 
\(Y\), and if \(\Phi \in \GKK^\Grd_j(Y, W)\) is another, then 
we write 
\( \Psi \otimes_Y \Phi\in \GKK^\Grd_{i+j}(X, W)\)
for their composition. 

We do not 
describe the general intersection product here, but will focus instead on 
the transversality method of \cite{Connes-Skandalis:Longitudinal}. 

Recall that two smooth \(\Grd\)-maps \(f_1\colon \Midd_1 \to  Y\) 
and \(b_2\colon \Midd_2 \to Y\) 
are 
\emph{transverse} if for every \( (p_1, p_2) \in \Midd_1\times \Midd_2\) 
 such that \(f_1(p_1)=b_2(p_2)\), we have \(D_{p_1}f_1(\Tvert_{p_1}\Midd_1) + 
 D_{p_2}b_2 (\Tvert \Midd_2) = \Tvert_{f_1(p_1) } (X)\). Transversality 
 ensures that 
 the space \[ M_1\times_X\Midd_2 \defeq \{(p_1, p_2) \in \Midd_1\times \Midd_2
 \mid  f_1(p_1)=b_2(p_2)\} \] has the structure of a smooth 
 \(\Grd\)-manifold.  

\begin{theorem}
  \label{thm:compose_transversely}
  Let \(\NM_1 = (\Midd_1,\mapl_1,\mapr_1,\Kclass_1)\) and
  \(\NM_2 = (\Midd_2,\mapl_2,\mapr_2,\Kclass_2)\) be smooth
  correspondences from~\(\Source\) to~\(\Target\) and
  from~\(\Target\) to~\(\Third\), respectively.  Assume that
  both \(\Midd_1\) and~\(\Midd_2\) admit smooth normally non-singular map s
  to~\(\Base\) (see \S \ref{subsec:normal_maps}),
   so that we lose nothing if we view \(\mapr_1\)
  and~\(\mapr_2\) as \(\K\)-oriented smooth maps (see 
  Theorem \ref{thm:normal_map_unique}).
  
    Assume also that \(\mapr_1\) and~\(\mapl_2\) are transverse, so that 
  \(\Midd_1\times_\Target\Midd_2\) is a smooth
  \(\Grd\)\nb-manifold; it has a smooth normally non-singular map  to~\(\Base\) as
  well, and the intersection product of \(\Phi_1\) and \(\Phi_2\) is 
  the class of the correspondence
  \[
 \bigl(\Midd_1\times_\Target\Midd_2, \mapl_1\circ\pi_1,
  \mapr_2\circ\pi_2, \pi_1^*(\Kclass_1)\cdot
  \pi_2^*(\Kclass_2)\bigr),
  \]
  where \(\pi_j\colon \Midd_1\times_\Target \Midd_2\to\Midd_j\)
  for \(j=1,2\) are the canonical projections.
\end{theorem}

In the non-equivariant situation, any two smooth maps can be perturbed 
to be transverse, and in \cite{Connes-Skandalis:Longitudinal} this is 
shown to give rise to a bordism of correspondences. As a result, one can 
compose bordism classes of correspondence by the recipe described in 
Theorem \ref{thm:compose_transversely}. 

However, this fails in the equivariant situation because 
 pairs of smooth maps cannot in general be 
perturbed equivariantly to be transverse; this happsn in even some of the 
simplest situations. 

\begin{example}
\label{ex:not_transverse}
Let \(\mu\) be the non-trivial character of \(\Z/2\). 
The corresponding one-dimensional representation 
is denoted \(\C_\mu\). We regard this as an 
equivariant vector bundle \(\C_\mu\) over a point. 
Its total space is \(\total{\C_\mu}\). 
The equivariant vector bundle \(\C_\mu\) 
is equivariantly \(\K\)-oriented, since 
the \(\Z/2\)-action preserves the complex structure on \(\C_\mu\). 
We therefore obtain a smooth normal equivariant map 
\(\star \to \total{\C_\mu}\) of degree \(2\). But since
the origin is the only fixed-point of the \(\Z/2\)-action, this map cannot be 
perturbed to be transverse to itself. 

This means that we cannot compose, for example, the topological 
morphism \(x\in \GKK^{\Z/2}_2(\star, \total{\C_\mu})\) 
represented by the correspondence \(
\star \leftarrow \star \rightarrow \total{\C_\mu} \), and the topological 
morphism 
\(y \in \GKK^{\Z/2}_0(\total{\C_\mu}, \star)\)
 represented by \(\total{\C_\mu} \leftarrow \star 
\rightarrow \star \) in the order \(x\otimes_{\total{\C_\mu}} y\), using the 
transversality recipe of Theorem \ref{thm:compose_transversely}. 

However, we may apply Thom modification to the 
 correspondence 
\(\star \leftarrow \star \rightarrow \total{\C_\mu}\), using \(\C_\chi\). Modification 
replaces the middle space \(\star\) of \(x\) 
by \(\total{\C_\chi}\) 
and replaces the normally non-singular map  \(\star\to \C\) by 
the composition \(\total{\C_\chi} \to \star \to \total{\C_\mu}\).
 This latter is obviously isotopic 
to the identity map on \(\total{\C_\mu}\). Finally, one adds the 
Thom class \(\xi_{\C_\mu}\in \K^{2}_{\Z/2}(\total{\C_\mu})\). 
Therefore the morphism \(x\) is equivalent to the morphism represented by 
\(\star \leftarrow (\total{\C_\mu}, \xi_{\C_\mu})  \xrightarrow{\ID} \total{\C_\mu}\). The 
identity map is transverse to any other map, since it is a submersion. 
Composing the modified correspondence and the original representative of 
\(y\) using the transversality recipe 
yields the class of the degree \(2\) correspondence   
\[\star \leftarrow (\star, (\xi_{\C_\mu})_{|_{\star}}) \rightarrow \star ,\]
where \((\xi_{\C_\mu})_{|_{\star}} \) denotes the restriction of the 
Thom class to the point. This equals the difference 
\([\epsilon] - [\chi]\in \Rep (\Z/2)\) of 
the trivial and the non-trivial representation of \(\Z/2\). 
It is the \emph{Euler class} of the \(\K\)-oriented vector bundle \(\C_\chi\),
 \emph{e.g.} the restriction of the Thom class to the zero section. 

\end{example}

With the given architecture of equivariant correspondences, 
a similar process can be carried out when composing two 
arbitrary (smooth) correspondences. Let 
\( X\xleftarrow{\mapl} (M, \xi) \xrightarrow{\mapr} Y\) be such. 
Let \(f = (\VB, \Triv, \hat{f})\). The equivariant vector bundles 
\(\VB\) and \(\Triv\) are \(\K\)-oriented by assumption. Thom modification 
using \(\VB\) results in the correspondence 
\( X\xleftarrow{\mapl} 
(\total{\VB}, \xi_\VB \cdot \xi) \xrightarrow{\mapr\circ \proj{\VB}} Y\). 
An obvious smooth isotopy of normally non-singular map s replaces this by 
\( X\xleftarrow{\mapl} 
(\total{\VB}, \xi_\VB \cdot \xi) \xrightarrow{\proj{\Triv^Y}\circ \hat{f}} Y\).  
Since \( \proj{\Triv^Y}\circ \hat{f}\) is a submersion, it is transverse to 
any other smooth map to \(Y\). Hence this correspondence 
can be composed using the transversality recipe of Theorem 
\ref{thm:compose_transversely} with any other one (on the right). 

An analogous procedure can be used to define a composition rule 
for arbitrary (not necessarily smooth) correspondences. 
This rule is quite topological in flavour, of course, 
but is only defined up to isotopy and is less satisfying than 
the sharp formulas one gets in the presence of traversality, which 
of course only apply in the presence of smooth structures. 
We will only compute compositions in this setting in this survey.

\subsection{Further properties of topological \(\KK\)-theory}
We have said that \(\GKK^\Grd\) is a category. It is also additive, 
with the sum operation on correspondences defined by a disjoint 
union procedure. The other important property is the existence of 
external products. This means that there exists an external product 
map 
\[ \GKK^\Grd_i (X,Y) \times \GKK^\Grd_j(U,V) \to
 \GKK^\Grd_{i+j} (X\times_\Base U, Y\times_\Base V).\]
It leads to the structure on \(\GKK^\Grd\) of a symmetric monoidal category. 

Finally, there is a natural map \(\GKK^\Grd\to \KK^\Grd\). This is defined not 
using the pseudodifferential calculus, as in \cite{Connes-Skandalis:Longitudinal}, 
but by purely topological considerations. Indeed, by definition, a normal 
\(\K\)-oriented \(\Grd\)-map 
from \(X\) to \(Y\) factors, by definition, as a composite of a zero section of an 
equivariantly \(\K\)-oriented
 vector bundle, an equivariant open embedding, and the projection 
map for another equivariantly \(\K\)-oriented vector bundle. Zero sections and 
bundle projections yield elements of \(\KK^\Grd\) because of the Thom isomorphism 
of \cite{LeGall:KK_groupoid}. Open embeddings clearly determine elements 
morphisms in \(\KK^\Grd\) because they determine equivariant 
*-homomorphisms. 

The various naturality properties of the Thom isomorphism 
imply corresponding facts about the map \(\GKK^\Grd\to \KK^\Grd\). 
Other functorial properties of \(\GKK^\Grd\), \emph{e.g.} with respect 
to homomorphisms \(\Grd'\to \Grd\) of groupoids, are explained in 
detail in \cite{Emerson-Meyer:Geometric_KK}. 


\section{Topological duality and the topological Lefschetz map}
\label{sec:dualities_and_lef}
We have organized this survey around the goal of computing the 
Lefschetz map for smooth \(\Grd\)-manifolds. 
This problem is intertwined with that of computing equivariant 
\(\KK\)-groups topologically and we will solve both problems 
at once in this section. The first step is to 
describe a class of \(\Grd\)-spaces 
to which the general theory of duality described in 
\S \ref{sec:abstract_duals} applies.  Subject to the resulting 
constraints on \(X\) we 
will obtain a \emph{topological model} of the Lefschetz map and 
simultaneously a proof that the \(\GKK^\Grd\to \KK^\Grd\) to be  
is an isomorphism on both the domain and range of the Lefschetz map. 
This will complete our computation of \(\Lef\) for 
a fairly wide spectrum of \(\Grd\)-spaces \(X\).

\subsection{Normally non-singular \(\Grd\)-spaces}
\label{subsec:normal_finite__type}
\begin{definition}
A \emph{normally non-singular \(\Grd\)-space} \(X\) is a \(\Grd\)-space 
equipped with a \(\Grd\)-equivariant 
normally non-singular map  \((\Normal_X, \Triv, \hat{g})\) from \(X\) to \(\Base\).

We also require of a normally non-singular \(\Grd\)-space \(X\) that 
every \(\Grd\)-equivariant vector bundle on \(X\) is subtrivial. 
 \end{definition}
The vector bundle \(\Normal_X\) is called the \emph{stable 
normal bundle of \(X\)}. 
 
 We may assume without loss of generality that \(\Triv\) is equivariantly 
 \(\K\)-oriented. Since the zero bundle is always uniquely \(\K\)-oriented, we obtain 
 an equivariant \(\K\)-oriented normally non-singular map  \( (0_{\Normal_X}, \Triv, \hat{g})\) from
  \(\Normal_X\) to \(\Base\). Let 
  \(D\in \GKK^\Grd(\Normal_X, \Base)\) be the corresponding class.

  Since \(\GKK^{\Grd\ltimes X}\) has external products, we can define a map 
  \[
   \GKK^{\Grd\ltimes X} (X\times_\Base U, Y\times_\Base V) \to 
  \GKK^{\Grd\ltimes X} (\Normal_X \times_\Base U, \Normal_X\times_\Base V)\]
  If every \(\Grd\)-equivariant vector bundle over \(X\) is subtrivial, then 
  there is a forgetful functor \(\GKK^{\Grd\ltimes X} \to \GKK^\Grd\) and this 
  results in a map 
  \[\GKK^{\Grd\ltimes X} (\Normal_X \times_\Base U, \Normal_X\times_\Base V)
  \to \GKK^\Grd((\Normal_X \times_\Base U, \Normal_X\times_\Base V).\]
  Composing with the previous one yields a topological analogue 
  \[ \GKK^{\Grd\ltimes X}(X\times_\Base U, X\times_\Base V) \to 
  \GKK^\Grd(\Normal_X\times_\Base U, \Normal_X\times_\Base V)\]
  of the functor 
  denoted \(T_\Dual\) in the discussion of duality in 
  \S \ref{sec:abstract_duals}, 
  and, composing further 
  with the morphism \(D\) gives a topological analogue of the Kasparov duality map 
  \[\PD^*\colon \GKK^{\Grd\ltimes X}(X\times_\Base U, X\times_\Base V) \to
  \GKK^\Grd(\Normal_X\times_\Base U, V)\]
  of Theorem \ref{thm:first_duality}. However, it is of course \emph{not} the 
  case in general that every equivariant vector bundle over \(X\) is subtrivial, 
  which is why we have added this as a hypothesis. 

\begin{example}
\label{ex:normal_not_subtrivial}
Let \(X \) be the integers with the trivial action of the circle group 
\(\Grd \defeq \mathbb{T}\). Then there are (\emph{c.f.} 
Example \ref{ex:not_subtrivial}) equivariant vector bundles 
on \(X\) which are not subtrivial. 
This is despite the fact that  
\(X\) admits a normally non-singular map  to a point, since it smoothly embeds in the 
trivial representation of \(\Grd\) on \(\R\), with trivial normal bundle.

In any case, \(X\) is not normally non-singular. 

 \end{example}

  Any smooth 
  \(\Grd\)-manifold satisfying one of the hypotheses of Theorem 
  \ref{thm:normal_factorisation} is normally non-singular. 
  
  To define a topological analogue of the map denoted \(\PD\) we 
  need a class 
  \( \Theta \in \GKK^{\Grd\ltimes X}(X, X\times_\Base \Normal_X)\). 
  Combining composition with this class and the map 
  \[ \KK^\Grd(\Normal_X\times_\Base U, V) \to 
  \GKK^{\Grd\ltimes X} (X\times_\Base \Normal_X\times_\Base U, X\times_\Base V),\]
  which is defined in the topological category in the same way as in the analytic one, 
  we will obtain the required topological map 
  \[\PD\colon \KK^\Grd(\Normal_X\times_\Base U, V) \to
  \GKK^{\Grd\ltimes X} (X\times_\Base U, X\times_\Base V).\]
  We just describe 
  \(\Theta\) in a heuristic fashion. 
   Assume for simplicity that the base 
  \(\Base\) of \(\Grd\) is a point. 
  So \(\Triv\) is just a Euclidean space. 
  Choose a point \(x\in X\). 
  Using the zero section of \(\Normal_X\) and 
  the map \(\hat{g}\), we see \(x\) as 
  a point in the 
  open subset \(\total{\Normal_X}\) of \(\Triv\), and hence by 
  rescaling \(\Triv\) (fibrewise) into a 
  sufficiently small open ball around \(x\) we obtain an open embedding of 
  \(\Triv\) into \(\Normal_X\). Explicitly, we use an open embedding of the form
  \[\hat{\varrho}_x\colon \total{\Triv} \xrightarrow{\beta_x} B_\epsilon\bigl( \hat{g}(\zers{\Normal_X}(x))\bigr) 
  \subset \total{\Normal_X}\] where the first map is the re-scaling.
    This yields an obvious 
  normally non-singular map  
  \( (\Triv, 0_{\total{\Normal_X}}, \hat{\varrho}_x)\) from a point to \(\total{\Normal_X}\)
   It is easily 
  checked that this can be carried out continuously 
  in the parameter \(x\in X\), and we obtain an \(X\)-equivariant normal 
  \(\K\)-oriented map \(\delta\) from \(X\) to \(X\times_\Base \total{\Normal_X}\) 
  with trace the graph of the zero-section \(\Xi \colon X \to X\times_\Base \total{\Normal_X} 
  , \; \Xi (x) \defeq (x, (x,0))\). 
  This yields an 
  element of \(\GKK^X(X, X\times \total{\Normal_X})\). The same can be checked 
  to work equivariantly, and a fibrewise version works for groupoids with nontrivial 
  base. We let 
  \(\Theta \in \GKK^{\Grd\ltimes X}(X, X\times_\Base \total{\Normal_X})\) be the 
  corresponding class. 
   
 \begin{theorem}
 \label{thm:normal_duality}
 In the notation above, let \(X\) be a normal \(\Grd\)-space of finite type, let 
 \(\Normal_X\) be the stable normal bundle and 
 \(D \in \KK^\Grd(\Normal_X, \Base)\) and \(\Theta 
 \in \KK^{\Grd\ltimes X} ( X, X\times_\Base \Normal_X)\) the classes 
 constructed above.
 
  Then  \( ( \Normal_X, D, \Theta)\) is a 
 Kasparov dual for \(X\) in \(\GKK^\Grd\). , and the maps 
 \(\PD\) and \(\PD^*\) are isomorphisms. 
 \end{theorem}
 
 The same formal computations as in \S \ref{sec:abstract_duals} then imply that
 the maps 
 \(\PD\) and \(\PD^*\) are isomorphisms. 
 
 Again, the stronger result is proved in \cite{Emerson-Meyer:Geometric_KK} 
 that one gets a symmetric 
 Kasparov dual; this, remember this is designed to give, as well, 
  an isomorphism of the form
 \[ \GKK_\Grd^* ( X\times_\Base U, V) \cong 
 \GKK_{\Grd\ltimes X}^* (X\times_\Base  U, \Normal_X\times_\Base V)\]
for any pair of \(\Grd\)-spaces \(U\) and \(V\). 
We do not give the details. As a consequence one deduces the following 
theorem by using duality to reduce from the bivariant to the 
monovariant case. 

\begin{theorem}
  \label{thm:top_complete}
  Let \(X\) be a normal \(\Grd\)-space of finite type, and \(\Target\) be an 
  arbitrary \(\Grd\)-space. Then 
  the natural transformation 
  \(\GKK^*_\Grd(\Source,\Target) \to
  \KK^\Grd_{*} \bigl(\CONT_0(\Source), \CONT_0(\Target)\bigr)\)
  is invertible. 
\end{theorem}

\begin{example}
Consider the space \(X \defeq \Z\) with the trivial action of 
\(\Grd\defeq \mathbb{T}\). This space is not normally non-singular, 
though it is a smooth \(\Grd\) manifold admitting a normally non-singular 
map to a point. The map 
\begin{equation}
\label{eq:map_not_an_isomorphism}
 \GKK^\Grd(X, \star) \to \KK^\Grd(\CONT_0(X), \C)\end{equation}
is \emph{not} an isomorphism in this case. By duality, 
the elements of \(\KK^\Grd(\CONT_0(X), \C)\) are parameterised by 
\(\Grd\)-equivariant complex vector bundles on \(X\). One can check that 
the elements of \(\GKK^\Grd(X, \star)\) are by contrast parameterised by 
\(\Grd\)-equivariant complex vector bundles which only involve 
a finite number of representations of \(\Grd\).  In other words, 
\eqref{eq:map_not_an_isomorphism} is equivalent to the embedding 
\[ \oplus_{n \in \Z} \Rep (\mathbb{T}) \to \prod_{n\in \Z} \Rep (\mathbb{T})\]
of the direct sum into the direct product of representation rings. 

If we dropped the subtriviality requirement on vector bundles that 
we imposed on cycles for \(\GKK\) then \eqref{eq:map_not_an_isomorphism}
would be an isomorphism, but then we would not be able to define the 
intersection product of correspondences in general. 

\end{example}

\begin{remark}
\label{rem:smooth_vs_non_smooth_versions}
Restricting to smooth correspondences and smooth 
equivalence classes of correspondences yields a parallel 
`smooth' theory. If \(X\) is a smooth normally non-singular \(\Grd\)-space 
and \(Y\) is smooth,  
then it follows that the smooth and non-smooth versions of
\(\GKK_\Grd^* (X,Y)\) agree. 

\end{remark}

We are now in a position to solve the problem we have been 
working towards: a topological computation of the Lefschetz map for a 
normal \(\Grd\)-space of finite type. 
Firstly, we define the \emph{topological Lefschetz map} 
\[ \Lef\colon \GKK^{\Grd\ltimes X}_*(X\times_\Base X, X) \to \GKK_\Grd^*(X, \Base)\]
for any normal \(\Grd\)-space of finite type as in Definition \ref{def:Lefschetz_map}, 
using the topological Kasparov dual constructed above out of the stable normal bundle. 
Since a topological dual maps to an analytic dual, the diagram 
that the diagram 
\begin{equation}
 \label{eq:two_lefs}
  \xymatrix{ 
  \GKK^{\Grd \ltimes \Tot}_*(\Tot\times_\Base \Tot  ,\Tot) 
  \ar[d] \ar[r]^{\Lef} & \GKK_\Grd^*(\Tot,\Base)\ar[d] \\
  \KK^{\Grd\ltimes \Tot}_*\bigl( \CONT_0(\Tot\times_\Base \Tot) ,
  \CONT_0(\Tot)\bigr) \ar[r]^{\;\;\;\;\;\;\;\Lef} & \KK_*^\Grd\bigl(\CONT_0(\Tot),\CONT_0(\Base)\bigr) 
   }.\end{equation}
 commutes. The finite-type hypothesis implies (Theorem \ref{thm:top_complete}) 
 that the vertical maps are both isomorphisms. Therefore we have 
 obtained a complete description of \(\Lef\) in purely topological terms. We 
 now proceed to describe it explicitly for smooth \(\Grd\)-manifolds in terms of 
 transversality. 
 
 \subsection{Explicit computation of \(\Lef\)}
Let \(X\) be a smooth \(\Grd\)-manifold satisfying 
 one of the hypotheses \(\bf{A}\) or \(\bf{B}\) of Theorem 
 \ref{thm:normal_factorisation}. From that Theorem, and by 
 Remark \ref{rem:smooth_vs_non_smooth_versions}, the 
 map \(X\to \Base\) is the trace of an essentially unique smooth 
 normally non-singular \(\Grd\)-map. The normal bundle \(\Normal_X\) is a 
 smooth \(\Grd\)-equivariant vector bundle, and 
 \(\hat{g}\colon \Normal_X \to \Triv\) is a smooth \(\Grd\)-equivariant 
 embedding.

 Moreover, the general morphism in \(\GKK^{\Grd\ltimes X}(X\times_\Base X, X)\)
 is represented by a smooth, \(\Grd\ltimes \Tot\)-equivariant correspondence
 \[ \Psi \defeq (M, \mapl, \mapr, \xi), \; \text{or in diagram form} \;
 X\times_\Base X\xleftarrow{\mapl} (M, \xi) \xrightarrow{\mapr} X\]
  from \(\Tot\times_\Base \Tot\) 
 to \(\Tot\). Let \(\Psi\) denote its class. 
 Recall that the \(X\)-structure on \(X\times_\Base X\) is on the 
 \emph{first} coordinate. 
 
 \begin{remark}
 The map \(\mapr\colon M \to X\) embedded in the correspondence 
 \(\Psi\) is assumed a 
 smooth \(\Grd\ltimes X\)-equivariant normally non-singular map . This presupposes the structure on \(M\) of a 
 smooth \emph{\(\Grd\ltimes X\)-manifold}. Note that this is a 
 stronger condition than 
 being a smooth \(\Grd\)-manifold: it entails a bundle structure
 on \(M\) with smooth fibres and is equivalent to requiring that 
 the smooth normally non-singular \(f\) is a submersion. 
 \end{remark}
 
 Following the definition of the Lefschetz map in
  \S \ref{subsec:the_lefschetz_map} 
 we next apply the functor 
 \(T_{\Normal_X}\colon \GKK^\Grd\to \KK^\Grd\) which sends a 
 \(\Grd\ltimes X\)-space, \emph{i.e.} a \(\Grd\)-space 
 \(W\) over \(X\), to the \(\Grd\)-space \(W\times_X \Normal_X\). The latter is
 the same as the pullback to \(W\) of the vector bundle \(\Normal_X\) using the 
 map \(W\to X\); recall that this functor is well-defined provided that 
 every \(\Grd\)-equivariant vector bundle over \(X\) is subtrivial.

 The functor \(T_{\Normal_X}\) maps a \(\Grd\ltimes X\)-equivariant
 map from \(W\) to \(V\) to the \(\Grd\)-equivariant 
 map \(W\times_X\Normal_X\to V\times_X\Normal_X\) given by the 
 obvious formula. 
 Since \(X\times_\Base X \times_X \Normal_X \cong X\times_\Base \Normal_X\) via 
 the map which forgets the first coordinate, applying the map \(T_{\Normal_X}\) to 
 \(\Psi\) yields 
 the \(\Grd\)-equivariant correspondence 
 \[ X\times_\Base \Normal_X\xleftarrow{\bar{b}}(M\times_X\Normal_X, \xi) \xrightarrow{\bar{\mapr}}\Normal_X,\]
 where \(\bar{b} (m, x,\xi) \defeq \bigl( b'(m), x,\xi)\), where \(b'\defeq \pr_2\circ b\) is the 
 second coordinate value of \(b\colon M\to X\times_\Base X\) (since the first 
 coordinate value is just the anchor map \(\anchor_M^{\Grd\ltimes X} \colon M \to X\) the 
 map \(\mapl\) is determined by \(\mapl'\).) 
 
 We then compose with the class of the smooth \(\K\)-oriented normally non-singular map  
 \(\Normal_X \to \Base\), using transversality (see 
 Theorem \ref{thm:compose_transversely}) in the category of 
 smooth, \(\Grd\)-equivariant correspondences. In 
 the notation of the theorem 
 \(M_1\defeq M\times_X \Normal_X\) and \(M_2 = \Normal_X\), and 
 \(f_1 = f\times_{X} \ID\colon
 M\times_X\Normal_X \to X\times_X \Normal_X = \Normal_X\), and 
 \(b_2\defeq \ID \colon \Normal_X \to \Normal_X\).  Since 
 \(\Normal_X \to \Base\) and \(f_1\) are smooth and normal 
 by assumption and \(f_1\) and \(b_2\) are obviously 
 transverse since \(b_2\) is a submersion, the composite using 
 Theorem \ref{thm:compose_transversely} is 
 the class of the smooth \(\Grd\)-equivariant correspondence 
 \begin{equation}
 \label{eq:next}
  X\times_\Base 
 \Normal_X\xleftarrow{\bar{b}}(M\times_X\Normal_X, \xi) \rightarrow \Base,
 \end{equation}
 and the only remaining step is to compose with the dual class, 
 which is the most interesting one from our point of view. 

  Recall that a bar denotes forgetting the \(X\)-structure on a \(\Grd\ltimes X\)-equivariant 
  morphism. The class \(\overline{\Theta}\) that we want to compose with 
  is that of the \(\Grd\)-equivariant correspondence 
 \[ X\xleftarrow{\ID} X\xrightarrow{\delta} X\times_Z \Normal_X\]
 -- which can be assumed smooth. Since the trace of 
 \(\delta\) is the smooth section \(\Xi\) where \(\Xi (x) = (x, (x,0))\),  one checks that 
 the transversality condition needed to compose \(\overline{\Theta}\) with 
 \eqref{eq:next} is that 
 \(\Xi\) is transverse to the map \(\bar{b}\). This can be checked to be
 the case if, for all \(m\in M\) for which \(\anchor_G^{\Grd\ltimes X} (m) = b'(m)\), 
  the linear map 
 \begin{equation}
 \label{eq:transverse_for_lef}
 \Tvert_mM \to \Tvert_{\anchor_M(m)}X, \hspace{1cm} \zeta \mapsto D_mb'(\zeta)-D_m\anchor_M^{\Grd\ltimes X}(\zeta)
 \end{equation}
 is non-singular. This thus implies that the 
 \emph{coincidence space} 
\[ \Fixed_\Psi \defeq \{m\in M \mid \rho_M^{\Grd\ltimes X}(m) = b'(m)\},\]
is a \emph{smooth, equivariantly \(\K\)-oriented 
\(\Grd\)-manifold}; in fact, more, by Theorem \ref{thm:compose_transversely}
it implies that the projection \(\Fixed_\Psi \to \Base\) is a smooth, 
\(\K\)-oriented normally non-singular map . Finally, it is easily checked that the restriction of 
\(\xi\) to \(\Fixed_\Psi\) has compact vertical support with respect to the map  
\((\anchor_M)_{_{|_{\Fixed_\Psi}}}^{\Grd\ltimes X}\colon \Fixed_\Psi  \to X\), (this is already 
implied by Theorem \ref{thm:compose_transversely}), 
so we get a \(\Grd\)-equivariant 
correspondence from \(X\) to \(\Base\): 
\[ (\Fixed_\Psi, (\anchor_M^{\Grd\ltimes X})_{|_{\Fixed_\Psi}}, \anchor_{\Fixed_\Psi}^\Grd, \xi),\]
in the usual notation for anchor maps. We call this the \emph{coincidence cycle} of \(\Psi\).

\begin{theorem}
\label{thm:main_lef}
Let \(\Psi  \in \GKK^{\Grd\ltimes X} _*(X\times_X, X)\); let 
\[\Lef\colon \GKK_\Grd^*(X\times_ZX, X) \to \GKK_\Grd^*(X, \Base)\]
be the Lefschetz map in topological equivariant \(\KK\)-theory. 
Then the topological Lefschetz invariant of the class of a correspondence 
\(\Psi\) in general position in the sense explained above, 
is the class of the coincidence cycle 
of \(\Psi\), 
\[\Lef \bigl( [(M, \mapl, \mapr, \xi)]\bigr) = 
[(\Fixed_\Psi, (\anchor_M)^{\Grd\ltimes X}_{|_{\Fixed_\Psi}}, \anchor_{\Fixed_\Psi}^\Grd, \xi )]
\in \GKK_\Grd^*(X, \Base).\]
Similar statements follow in analytic \(\KK\). 
\end{theorem}

Namely, the Lefschetz invariant
of the \(\KK^{\Grd\ltimes X}(X\times_\Base X, X)\)-morphism \(\KK^{\Grd\ltimes X} (\Psi)\) 
determined by \(\Psi\) is the class \(\KK^\Grd ( \Lef (\Psi))\). 
Furthermore, this class is the 
pushforward under the map
\((\anchor_M)^{\Grd\ltimes X}_{|_{\Fixed_\Psi}}\colon \Fixed_\Psi \to X\) of the 
class of the Dirac operator on the \(\K\)-oriented 
coincidence manifold \(\Fixed_\Psi\), twisted by \(\xi\) (by an appeal 
to the Index Theorem.) 

We leave it to the reader to compute the Lefschetz invariant of \(\Psi\) in 
the situation where the transversality condition \eqref{eq:transverse_for_lef} fails; 
in this case it becomes necessary to modify \(\overline{\Theta}_{\textup{top}}\) 
using Thom modification and an isotopy as in Example \ref{ex:not_transverse}. 

Our computation of 
the Lefschetz map for smooth normal \(\Grd\)-manifolds of finite type is now complete, 
in view of \eqref{eq:two_lefs} and Theorem \ref{thm:top_complete} in combination 
with Remark \ref{rem:smooth_vs_non_smooth_versions}.

\subsection{Lefschetz invariants of self-morphisms of \(X\)}
Let \(X\) be a smooth normal \(\Grd\)-manifold of finite type. 
We now consider the composition 
\[ \GKK^\Grd(X,X) \to \GKK^{\Grd\ltimes X}(X\times_\Base X, X)
\xrightarrow{\Lef} \KK^\Grd(X, \Base),\]
where the first map is the composition of the canonical inflation map 
\[p_X^*\colon \GKK^\Grd(X, X) \to \GKK^{\Grd\ltimes X} (X\times_\Base X, X\times_\Base X)\]
and the map 
\[ \GKK^{\Grd\ltimes X} (X\times_\Base X, X\times_\Base X)
\to \GKK^{\Grd\ltimes X}(X\times_\Base X, X)\]
of composition with the diagonal restriction class 
\( \Delta_X\in \KK^{\Grd\ltimes X}( X\times_\Base X, X).\) Recall that the latter 
is the class of the \(\Grd\ltimes X\)-equivariant correspondence 
\( X\times_\Base X \xleftarrow{\delta_X}X \xrightarrow{\ID}X\) where \(\delta_X\) is 
the diagonal embedding.  
Let \(\Psi = (M, \mapl, \mapr, \xi)\) be a smooth, 
\(\Grd\)-equivariant correspondence 
from \(X\) to \(X\). Of course this implies that 
there is a smooth normally non-singular map  \(M \to \Base\), by composing 
\(\mapr\colon M \to X\) and 
\(\anchor_X\colon X \to \Base\). The inflation map replaces \(\Psi\)  
by the \(\Grd\ltimes X\)-equivariant correspondence 
\[X\times_\Base X \xleftarrow{\ID_X\times_\Base \mapl} (X\times_\Base M, \xi) 
\xrightarrow{\ID\times_\Base \mapr} X\times_\Base X.\]
The \(X\)-structures are all on the first variable. In order to compose 
(on the right) with the diagonal restriction class using transversality, we 
first easily check that Theorem \ref{thm:compose_transversely})
 applies, and deduce that we require that 
the smooth maps 
\(\ID_X\times_\Base \mapr\colon X\times_X M \to X\times_X X\) and 
\(\Delta_X\colon X \to X\times_\Base X\) are transverse in the sense of 
Theorem \ref{thm:compose_transversely}, in the category of 
\(\Grd\ltimes X\)-equivariant smooth maps. \emph{They are transverse if and only 
the smooth \(\Grd\)-map \(\mapr\colon M \to X\) is a submersion.} If 
this condition is met, then \(\mapr\colon M \to X\) gives \(M\) not just the 
structure of a smooth \(\Grd\)-manifold, but the structure of a smooth 
\(\Grd\ltimes X\)-manifold, \emph{i.e.} a bundle of smooth 
manifolds over \(X\) with morphisms in \(\Grd\) acting by diffeomorphisms 
between the fibres.  Composing with \(\Delta_X\) using transversality 
then yields the \(\Grd\ltimes X\)-equivariant correspondence 
\( X\times_\Base X \xleftarrow{(\mapl, \mapr)} (M, \xi) 
\xrightarrow{\mapr} X\) where \( (\mapl, \mapr) (m)
 \defeq ( \mapl (m), \mapr (m))\).  Finally, we apply 
 Theorem \ref{thm:main_lef} to 
obtain the following. 

\begin{theorem}
\label{thm:lef_self_morphisms} 
Let \(X\) be a smooth normally non-singular \(\Grd\)-manifold,
let \(\Psi = (M, \mapl, \mapr, \xi)\) be a smooth, 
\(\Grd\)-equivariant correspondence 
from \(X\) to \(X\) such that \(f\) is a fibrewise submersion and 
such that for every \(m\in \Midd\) such that \(\mapl (m) = 
\mapr (m)\) the linear map \(\Tvert_mM \to \Tvert_{\mapl(m)}X\), 
\[ \zeta \mapsto D_m\mapl (\zeta) - D_m \mapr (\zeta)\]
is non-singular. 
Then 
the Lefschetz invariant of \(\Psi\) is the class of the smooth, 
\(\Grd\)-equivariant correspondence 
\[ X\xleftarrow{\mapl} (\Fixed_\Psi' ,\xi_{|_{\Fixed_\Psi'}})\rightarrow \Base\]
from \(X\) to \(\Base\),
where \(\Fixed_\Psi' \defeq \{m \in M \mid \mapl (m) = \mapr (m)\}\) with its 
induced \(\K\)-orientation. 
\end{theorem}

\begin{example}
\label{ex:lef_for_maps}
If \(\mapl\colon X \to X\) is a smooth \(\Grd\)-equivariant map, then 
\(X\xleftarrow{\mapl} (X, \xi) \xrightarrow{\ID} X\) is a smooth,  
zero-dimensional correspondence from \(X\) to \(X\). Since 
\(\ID\colon X \to X\) is obviously a submersion, the transversality condition 
amounts to saying that for every \(x \in X\) which is fixed by \(\mapl\), 
the linear map \(\ID - D_x\mapl\colon \Tvert_xX \to \Tvert_x X\) is 
non-singular. This is the classical condition. The coincidence space 
is of course the fixed-point set of \(b\), suitably \(\K\)-oriented. This 
means a sign is attached to each point, which can be checked to 
agree with the usual assignment. Thus the Lefschetz invariant 
is the algebraic fixed-point set. 
\end{example}

\subsection{Homological invariants for correspondences}
\label{subsec:homological_lef}
In this section we describe, in the non-equivariant case, the 
pairing  between the 
index of the Lefschetz invariant of a Kasparov self-morphism 
of \(X\), and  a
\(\K\)-theory class, in terms of graded traces on \(\K\)-theory. 
In combination with the theory of correspondences this yields 
a generalisation of the classical Lefschetz fixed-point formula. 

Recall that the classical Lefschetz fixed-point theorem describes the 
algebraic number of fixed-points of a map satisfying a transversality 
condition, to the graded trace of the induced map on homology, a 
`global' invariant. The latter of course only makes sense 
when the homology groups of \(X\) have finite rank.  
Duality and the Universal Coefficient Theorem taken together 
\emph{imply this}, so in particular it is the case if
 \(X\) is a compact manifold. All the same remarks hold of course for 
 \(\K\)-theory as well. 
 
 \emph{The reader may take all \(\K\)-theory groups to 
 be tensored by \(\Q\) in the following}. Recall that the \(\K\)-theory 
 \(\K^*(X)\) is a graded ring. This ring structure is important for what 
 follows.  Let \(L_x\colon \K^*(X)\to \K^*(X)\) 
 denote the \emph{additive} group homomorphism of \(\K^*(X)\) 
  by multiplication by \(x \in \K^* (X)\). 
  For \(f\in \KK_*(\CONT(X), \CONT(X))\), let
  \(f_*\) denote the endomorphism of \(\K\)-theory induced by \(f\). 
   We are interested in the linear transformation 
  \[ L(f) \colon \K^* (X) \to \Z, \hspace{1cm} 
  L(f)x \defeq \textup{trace}_s (f_*\circ L_x).\]
  We call \(L(f)\) 
  the \emph{Lefschetz operator}  of \(f\). It is a globally defined 
  object, generalizing the classical 
  \emph{Lefschetz number} \(l(f) \defeq \mathrm{trace}_s (f_*)\) of \(f\) 
  in the sense that evaluating \(L(f)\) at the unit \([1]\in \K^0 (X)\) 
  recovers the Lefschetz number:
  \[ L(f) ( [1] ) = l(f).\]
  The Lefschetz operator contains more information; if 
  for example if \(f\) is an \emph{odd} morphism then 
  \(\ell (f) = 0\) 
  but \(L(f)\not= 0\) except in special cases.

 \begin{theorem}
 \label{thm:homological_lef}
  Let \(X\) be a compact space admitting an abstract dual.  
  Let \(f\in \KK_*\bigl( \CONT(X), \CONT(X))\) be a Kasparov morphism. 
 Then 
 \begin{equation}
 \label{eq:homological_lef}
  L(f)\xi = \langle \xi, \Lef (f)\rangle  
  \end{equation}
  holds for every \(\xi\in \K^*(X)\).
 
 \end{theorem}
 
 By the characteristic-class formulation \cite{Atiyah-Singer:Index_III}
 of the Atiyah-Singer 
 Index theorem, for each compact, smooth \(\K\)-oriented 
  manifold \(M\) there exists a cohomology class \(\Index (M)\) on 
  \(M\) which we call the \emph{orientation character of \(M\)} such that 
  \[ \ind (D\cdot \xi) = \int_M \Index (M) \textup{ch}(\xi),\]
 where \(D\cdot \xi\) denotes the class of the Dirac operator on \(M\)
  twisted by \(\xi\). Of course we choose representative differential 
  forms for \(\Index (M)\) and \(\textup{ch}(\xi)\) 
  both can be computed more or 
  less explicitly from Chern Weil theory. 
 In combination with Theorem \ref{thm:homological_lef} we 
 obtain a local formula for the Lefschetz operator of a morphism 
 represented by a correspondence, as follows.  

 \begin{corollary}
 Let \( X\xleftarrow{\mapl} (M, \xi) \xrightarrow{\mapr} X\) 
 be a smooth self-correspondence of \(X\), assume it is of degree 
 \(d  \defeq \dim (M) -\dim (X) +\dim (\xi)\) and denote its class by 
 \(\Psi\). 
 
 Assume that the transversality assumptions in 
 Theorem \ref{thm:lef_self_morphisms} are met, so that 
 the coincidence space
  \[\Fixed_\Psi' \defeq \{m \in M \mid \mapl (m) = \mapr (m)\} \] 
  has the structure of a smooth, \(\K\)-oriented manifold of 
 dimension \(d\). Then  
  \begin{equation}
  \label{eq:lef_characteristic_class}
  L(f) \eta =
  \int_{ \Fixed_\Psi } 
  \textup{ch} 
  \bigl( \xi_{|_{\Fixed_\Psi'}}\cdot (\mapl^*\eta)_{|_{\Fixed_\Psi'}}\bigr)
  \Index(D_{\Fixed_\Psi}),
 \end{equation}
 holds, where \(\Index (D_{\Fixed_\Psi'})\) is the 
 index character of \(\Fixed_\Psi'\), and 
 \(L(f) \eta \defeq \mathrm{trace}_s (f_*\circ L_x) \) is the 
 Lefschetz operator applied to \(\eta\). 
 \end{corollary}
 
 In the case where \(M = X\), \(\mapr = \ID\) and 
 \(\mapl \colon X \to X\) is a smooth equivariant map in 
 general position, the coincidence manifold \(\Fixed_\Psi'\) is a 
 finite set of points, and \eqref{eq:lef_characteristic_class} reduces to 
 the traditional Lefschetz fixed-point theorem: 
 \begin{equation}
 \label{eq:zero_dimensional_lef}
 \mathrm{trace}_s (\mapl_*) =  \sum_{x\in \mathrm{Fix} (\mapl)} 
 \sign \det  (\ID - D_xf)  .\end{equation}
 
Ralf Meyer and I are currently aiming at an equivariant version of the 
above, but the work is not yet complete. 
Let \(\Grd\) be compact group.
 Let \(X\) be a compact \(\Grd\)-space. Then 
the \(G\)-equivariant \(\K\)-theory \(\K^*_G(X)\)of \(X\) is a module 
over the representation ring  \(\Rep (G)\). The \emph{Hatori-Stallings trace}
\[ \textup{trace}_{\Rep (G)} \colon \K^*_G(X) \to \Rep (G)\]
is then defined under suitable conditions. Now assume that \(X\) admits an 
abstract dual, so that the Lefschetz map is defined. 
 If \(f\in \KK^*(\CONT_0(X), \CONT_0(X))\) is 
a morphism, then call \emph{the Lefschetz index} 
\(\textup{ind}_G\circ \Lef (f)\in \Rep (G)\) the pairing of 
the unit class \([1]\in \K^0_G(X)  = \KK^G(\C, \CONT_0(X))\) with 
\(\Lef (f) \in \KK^G(\CONT_0(X), \C)\). Our expectation is 
that the following result holds -- we state it as a conjecture 
since the proof is not complete at the time of writing.

\begin{conjecture}
\label{conj:equivariant_lefschetz}
If \(f\in \KK^G_*(\CONT_0(X), \CONT_0(X))\), then 
\(\textup{ind}_G\circ \Lef (f) = \textup{trace}_{\Rep (G)} (f_*),\)
where \(f_*\) denotes the action of \(f\) on \(\K^*_G(X)\) and 
\(\textup{trace}_{\Rep (G)}\) denotes the Hattori-Stallings trace. 
\end{conjecture}

This theorem can be of course combined with correspondences to 
achieve interesting local-global equalities of what seem to be 
rather subtle invariants. 

The Hattori-Stallings trace is defined for modules \(M\) over 
rings \(R\) which are \emph{finitely presented}, \emph{i.e.} have 
finite-length resolutions by finitely generated projective 
\(R\)-modules. Therefore we require this hypothesis on the 
equivariant \(\K\)-theory of \(X\) as a module over \(\Rep (G)\). 

This a reasonable assumption only for 
compact, connected Lie groups with 
torsion-free fundamental group \(G\).  For disconnected compact
Lie groups, \emph{i.e.} finite groups, 
the homological dimension of \(\Rep (G)\) 
typically has infinite homological dimension and we
do not know at the moment how to formulate an equivariant Lefschetz theorem for finite 
groups along the lines of \ref{conj:equivariant_lefschetz}.

\end{document}